\documentclass{siamltex}

\usepackage{amsfonts,amsmath,amssymb,graphicx}

\newcommand{\RR}{{\mathbb{R}}}
\newcommand{\CC}{{\mathbb{C}}}
\newcommand{\vect}{{\rm vec}}
\newcommand{\XX}{{\mathbf X}}

\newcommand{\xx}{{\mathbf x}}

\newcommand{\HH}{{\tt H}}
\newcommand{\TT}{{\tt T}}

\newtheorem{remark}[theorem]{Remark}

\newtheorem{example}[theorem]{Example}

\begin{document}
\bibliographystyle{plain}
\pagestyle{myheadings}
\markboth{V. Simoncini}{The Lyapunov matrix equation. Matrix analysis from a computational
perspective}

\title{The Lyapunov matrix equation. Matrix analysis from a computational 
perspective\thanks{Version of
January 10, 2015. This work is a contribution to the Seminar series ``Topics in Mathematics'', of
the PhD Program of the Mathematics Department, Universit\`a di Bologna.}}
\author{
V. Simoncini\thanks{Dipartimento di Matematica,
Universit\`a di Bologna,
Piazza di Porta San Donato  5, I-40127 Bologna, Italy
(valeria.simoncini@unibo.it).}
}
\maketitle
\begin{abstract}
Decay properties of the solution $X$ to the Lyapunov matrix equation
$AX + X A^\TT = D$ are investigated. Their exploitation in the
understanding of equation matrix properties, and in the
development of new numerical solution strategies 
when $D$ is not low rank but possibly sparse is also briefly discussed.
\end{abstract}

\begin{keywords}
sparsity pattern, banded matrices, Kronecker products, exponential decay.
\end{keywords}

\begin{AMS}
65F50, 15A09
\end{AMS}

\section{Introduction}
We are interested in the analysis of the linear matrix equation
\begin{eqnarray}\label{eqn:Lyap}
A X + X A^{\TT} = D , \qquad A, D\in\RR^{n\times n}, 
\end{eqnarray}
to be solved for $X \in\RR^{n\times n}$; here and in the following $A^\TT$ denotes
the transpose of the matrix $A$.
In particular, we focus on the decay and sparsity properties of the involved matrices that
can be exploited for computational purposes, or that can give insight
into the analysis of numerical solution methods.

Matrix equations have always had a key role in Control theory, because their
solution matrix carries information on the stability of the problem \cite{Antoulas.05},%
\cite{Benner2005a}.
More recently, linear matrix equations and their generalizations, linear ``tensor''
equations, have been shown to be an appropriate tool 
to represent the discretization of parameterized partial differential equations,
as they arise for instance in stochastic modeling; see, e.g.,
\cite{BabuskaetalSINUM.04},\cite{Khoromoskij.12},\cite{Furnival.Elman.Powell.10} 
and the discussion in \cite{Simoncini.survey13}. Their analysis and numerical
solution is therefore attracting considerable attention in the numerical and
engineering communities.

Using the Kronecker product, the matrix equation (\ref{eqn:Lyap})
can be rewritten
as the following standard (vector) linear system
\begin{eqnarray}\label{eqn:Kron}
{\cal A} \xx  = b, \qquad \mbox{with}\quad
\begin{array}{l}
{\cal A} =I_n\otimes A + A \otimes I_n \\
 x = \vect(X), \quad b=\vect(D), 
\end{array}
\end{eqnarray}
where  the Kronecker product of two matrices $X$ and $Y$ of size $n_x\times m_x$ and $n_y\times m_y$,
respectively, is defined as
 $$
 X \otimes Y = 
\begin{bmatrix}
x_{11}Y & x_{12} Y & \cdots & x_{1 m_A} Y \\
x_{21}Y & x_{22} Y & \cdots & x_{2 m_A} Y \\
 \vdots &          &        & \vdots      \\
x_{n_x 1}Y & a_{n_x 2} Y & \cdots & x_{n_x m_x} Y \\
\end{bmatrix} \in \CC^{n_x n_y \times m_x m_y};
 $$
the vec operator stacks the columns of a matrix $X=[x_1, \ldots, x_m] \in \CC^{n\times m}$
one after the other as
$$
{\rm vec}(X) = \begin{bmatrix} x_1 \\ \vdots \\ x_m\end{bmatrix} \in \CC^{nm} .
$$

From (\ref{eqn:Kron}) we can deduce that the system admits a solution for
any $b$ and this is
unique, if and only if the matrix $\cal A$  is nonsingular.
Using a standard result for the spectrum of Kronecker product matrices,
 this is equivalent to requiring that
${\rm spec}(A) \cap {\rm spec}(-A) = \emptyset$, where
${\rm spec}(A)$ denotes the set of eigenvalues of $A$ (see, e.g., \cite[Theorem~4.4.6]{Horn.Johnson.91}).

Though the Kronecker form (\ref{eqn:Kron}) of the problem is appealing, as a large
body of literature on linear systems can be exploited, this approach
dramatically increases the complexity of the computation, and also 
cannot preserve the intrinsic properties of the problem in practice; for instance,
if $D$ is symmetric, then the matrix $X$ is also symmetric. This property is
not preserved when solving (\ref{eqn:Kron}), unless the system is solved at
a very high accuracy.

The numerical solution of (\ref{eqn:Lyap})  is particularly challenging
when $n$ is large. Indeed, although $A$ and $D$ may be sparse and/or structured,
the solution matrix $X$ is usually dense. For $n=O(10^5)$ or higher, storing
the full matrix $X$ becomes prohibitive, and either sparse or low-rank
approximations are sought. Experimental evidence and theoretical results indicate 
that low rank approximations can indeed be sought after whenever the right-hand
side $D$ is itself a low rank matrix.  Major efforts in the past decades
have been devoted to determining approximations of the lowest possible rank,
given a fixed final accuracy. Interestingly, strategies that determine the
most accurate approximation for a given rank have also been investigated. 
We refer to \cite{Simoncini.survey13} for a recent
survey on recent computational strategies.

Sparsity and decay (quasi-sparsity) patterns of the solution matrix $X$ have 
been less investigated. We aim to explore recent findings in this direction:
while decay pattern results are proved by using the closed forms above and thus
they seem to be of merely theoretical interest, they
may be insightful in the development of computational procedures and in the
analysis of the expected solution properties.

Decay properties have been classically investigated for the entries
of the inverse of banded matrices. Let $A=(a_{i,j})$, $i,j=1\ldots, n$
 be a symmetric positive definite
banded matrix with bandwidth $\beta$, that is
$a_{i,j} = 0$ for $|i-j|>\beta$, and let $\kappa>1$ be the ratio between
its largest and smallest eigenvalues. Then Demko et al. \cite{SDWFMPWS84}
showed that
\begin{eqnarray}\label{eqn:Demko}
|A^{-1}|_{i,j} \le c_0 q^\frac{|i-j|}{\beta}
\end{eqnarray}
where $q = (\sqrt{\kappa}-1)/(\sqrt{\kappa}+1)<1$, and $c_0$ also depends
on the extreme eigenvalues of $A$. This bound shows an exponential off-diagonal
decay property as the inverse matrix entry moves away from the
main diagonal, and the decay depends on the bandwidth of $A$.
This property was generalized to functions of symmetric matrices in
\cite{Benzi.Golub.99} and more recently to more general matrices and
abstract settings in \cite{Benzi2007},\cite{Benzi.Boito.14}.

With completely different purposes (see \cite{CanutoSimonciniVeraniJSC.14}), 
in \cite{CanutoSimonciniVeraniLAA.14} it
was shown that also the Lyapunov operator
$$
{\cal L} \, : X \, \to \, AX + X A^\TT
$$
enjoys similar properties. More precisely, for $A$ symmetric and positive
definite, the entries of the 
{\it matrix} ${\cal L}^{-1}(e_i e_j^\TT)$ show a rapid, although not necessarily
exponential, decay pattern as one moves away from the element in position $(i,j)$.
Here we would like to linger over this property, exploring some of its
possible consequences both theoretically and computationally.
The main aim of this paper is to discuss some recent results, and highlight
directions of research that could improve our understanding of the
matrix properties of (\ref{eqn:Lyap}), and perhaps encourage the development of new 
computational strategies for its numerical solution when $D$ is not 
low rank but possibly sparse.
 
We will first generalize the decay pattern result to a wider class of matrices of interest
in applications. 
 Then we will 
show how to use this decay pattern to characterize the approximate solution matrix
of large Lyapunov equations
by projection-type methods. Finally, we will make some comments on the
solution pattern for sparse $D$.

All reported experiments and plots were performed using Matlab \cite{Matlab7}.

\section{Closed forms of the matrix solution}
The solution $X$ of (\ref{eqn:Lyap})
may be written in closed form in a number of different ways; here we report the
forms that have been used in the literature:

\begin{itemize}
\item[$(a)$] {\it Integral of resolvents.}
The following representation, due to Krein,
exploits spectral theory arguments: (see, e.g., \cite{Lancaster1970})
\begin{eqnarray}\label{eqn:integral}
X = \frac 1 {2\pi} \int_{-\infty}^{\infty}
(\omega \imath I -A)^{-1} D 
(\omega \imath I - A)^{-\HH} d\omega  ,
\end{eqnarray}
where $I$ is the $n\times n$ identity matrix
and $\Gamma_1$ is a  contour containing and sufficiently
close to the spectrum of $A$.

\item[$(b)$] {\it Integral of exponentials.}
Assume that the field of values\footnote{The field of values of $A\in\RR^{n\times n}$ 
is defined as $W(A) = \{ z \in\CC, \, z=x^\HH A x, x\in\CC^n, x^\HH x=1\}$.}
 of $A$ is all contained either in $\CC^-$ or
in $\CC^+$, excluding the imaginary axis.
This representation, due to Heinz,  \cite[Satz 5]{Heinz.55}, 
is tightly connected to the previous one,
\begin{eqnarray}\label{eqn:cauchy}
\XX =  \int_0^{\infty} e^{At} D 
e^{tA^\HH} dt , 
\end{eqnarray}
where $e^{At}$ is the matrix exponential of $At$.

\item[$(c)$] {\it Finite power sum.}
Assume $D=BB^\TT$ with $B\in\RR^{n\times s}$, $s\le n$, and let
$a_m$  of degree $m$ be the minimal polynomial  of $A$ with respect to
$B$, namely the smallest degree monic polynomial such that $a_m(A)B=0$.
 Then (\cite{Souza1981})
\begin{eqnarray}\label{eqn:finitesum}
\XX &=& 
\sum_{i=0}^{m-1}
\sum_{j=0}^{m-1}
\gamma_{ij} A^i C (A^\HH)^j \nonumber \\ 
&=&
[B, A B, \ldots, A^{m-1}B] (\gamma \otimes I) 
\begin{bmatrix} B^{\HH} \\ B^{\HH} A^{\HH} \\ \vdots \\ B^{\HH} (A^{\HH})^{m-1}\end{bmatrix} ,
\end{eqnarray}
where $\gamma$ is the solution of the Lyapunov equation with
coefficient matrices given by
the companion matrix of $a_m$, and right-hand side the matrix $e_1 e_1^\TT$, where
$e_1^\TT = [1, 0, \ldots, 0]$;
see also \cite{Lancaster.Lerer.Tismenetsky.84}.

\item[$(d)$] {\it Similarity transformations.} Strictly related to (c),
in addition this form assumes that $A$ can be
diagonalized, $U^{-1}AU = {\rm diag}(\lambda_1, \ldots, \lambda_n)$.
Let $\widetilde D = U^{-1} D U^{-\HH}$. Then
$$
X=U\widetilde X U^{-1}, \quad \mbox{with} \quad
\widetilde {x}_{ij} = \frac{ \widetilde D_{ij}}{\lambda_i+\mu_j}.
$$
\end{itemize}

These closed forms are not used for computational purposes, since stable
and efficient methods have been derived already in the 1970s \cite{Simoncini.survey13}.
Nonetheless, they can be used to describe other properties, such as
entry decay - as we are going to show - and numerical
low rank of the solution matrix \cite{Simoncini.survey13}.

\section{Decay pattern properties of the solution}
The sparsity pattern of the matrix $D$ influences the pattern of
the solution to (\ref{eqn:Lyap}), which can be qualitatively described
a-priori.
This pattern is typical of certain discretizations of
elliptic partial differential equations. For instance, when discretizing
the Poisson equation $-\Delta u = f$ on the unit square by means of 
finite differences or uniform low degree finite elements, 
the algebraic linear system ${\cal A} x = b$ is obtained, where
${\cal A} = A\otimes I + I\otimes A$ is called the stiffness
matrix, and $A$ is symmetric and
tridiagonal, $A={\rm tridiag}(-1,2,-1)$; see, e.g.,
 \cite[section 6.3.3]{Demmel.97}. Clearly, this problem
is equivalent to the Lyapunov equation (\ref{eqn:Lyap}) with
$x={\rm vec}(X)$ and $b={\rm vec}(D)$. The particular pattern we
are analyzing is derived whenever $f$ is, say, a single-point forcing term.
Analogously, the analysis of the decay pattern in $X$ 
may give insight into the understanding of the sparsity pattern
of the stiffness matrix $\cal A$, since each column $t$ of ${\cal A}^{-1}$,
${\cal A}^{-1} e_t$ is nothing but the solution to ${\cal A} x = e_t$
 \cite{CanutoSimonciniVeraniLAA.14}. Our analysis helps formalize the oscillating - while decaying -
pattern observed in the inverse stiffness matrix.

\begin{figure}[tb]
\centering
\includegraphics[width=1.35in,height=1.35in]{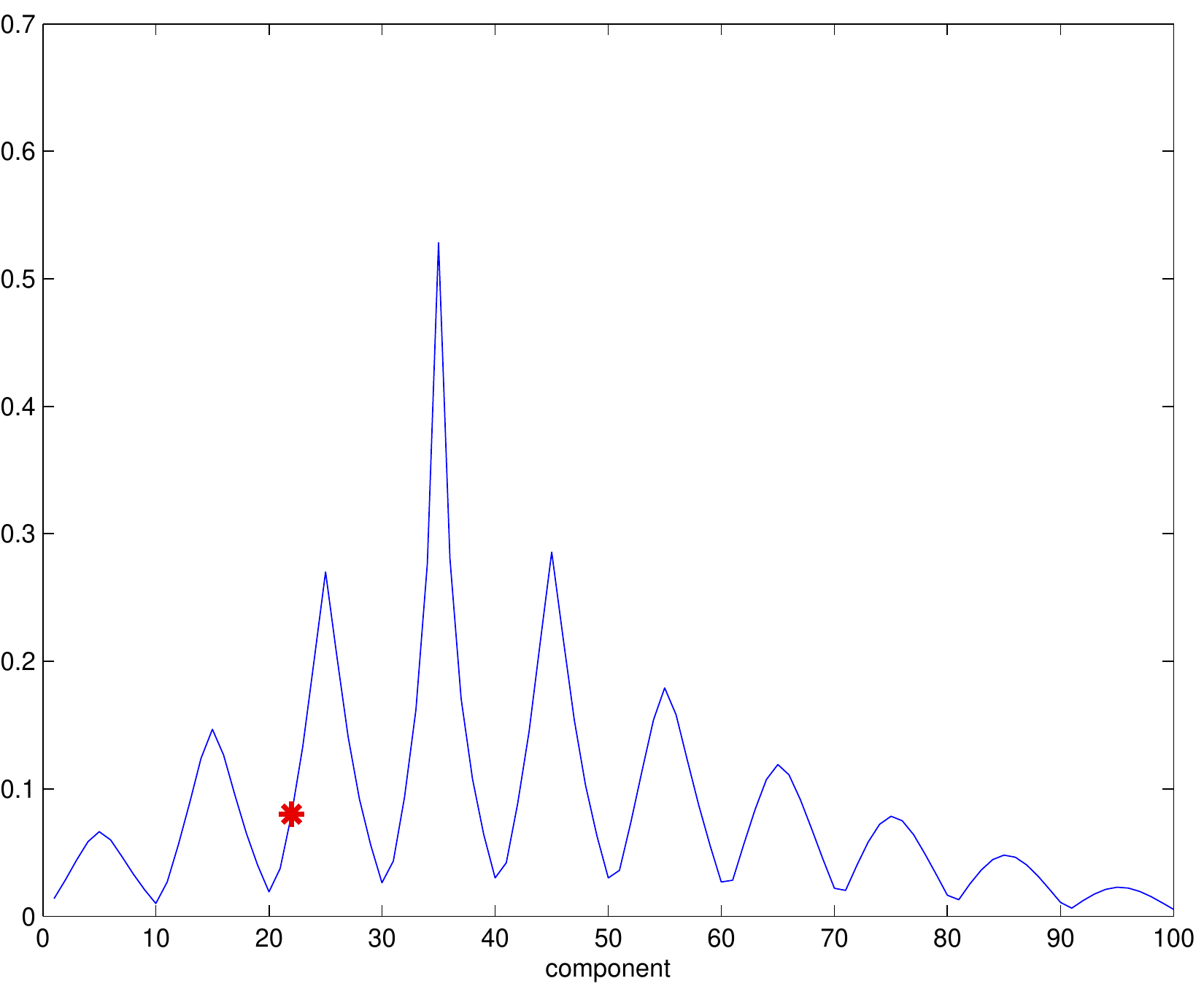}\hskip 0.5in
\includegraphics[width=1.35in,height=1.35in]{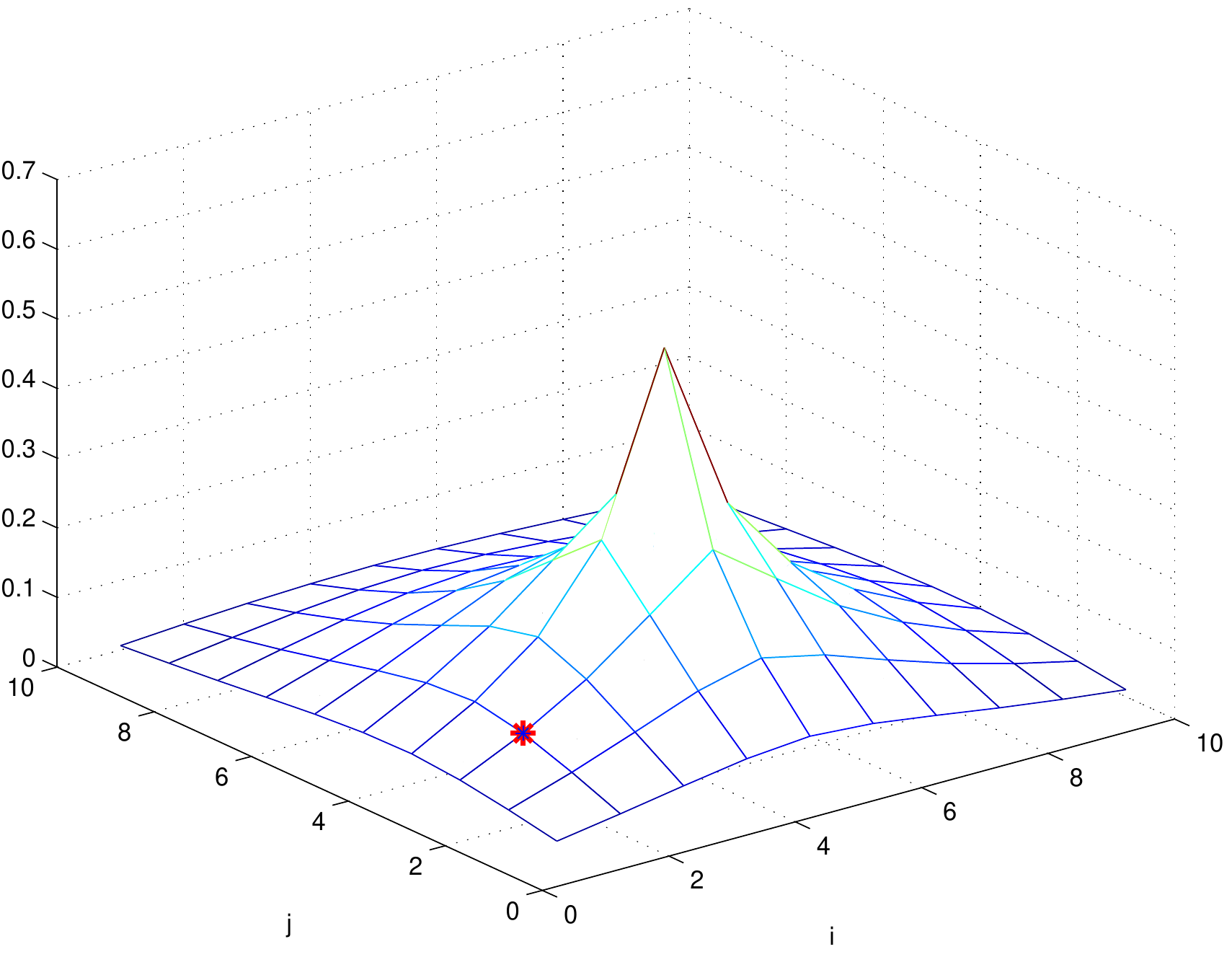} \\
\includegraphics[width=1.35in,height=1.35in]{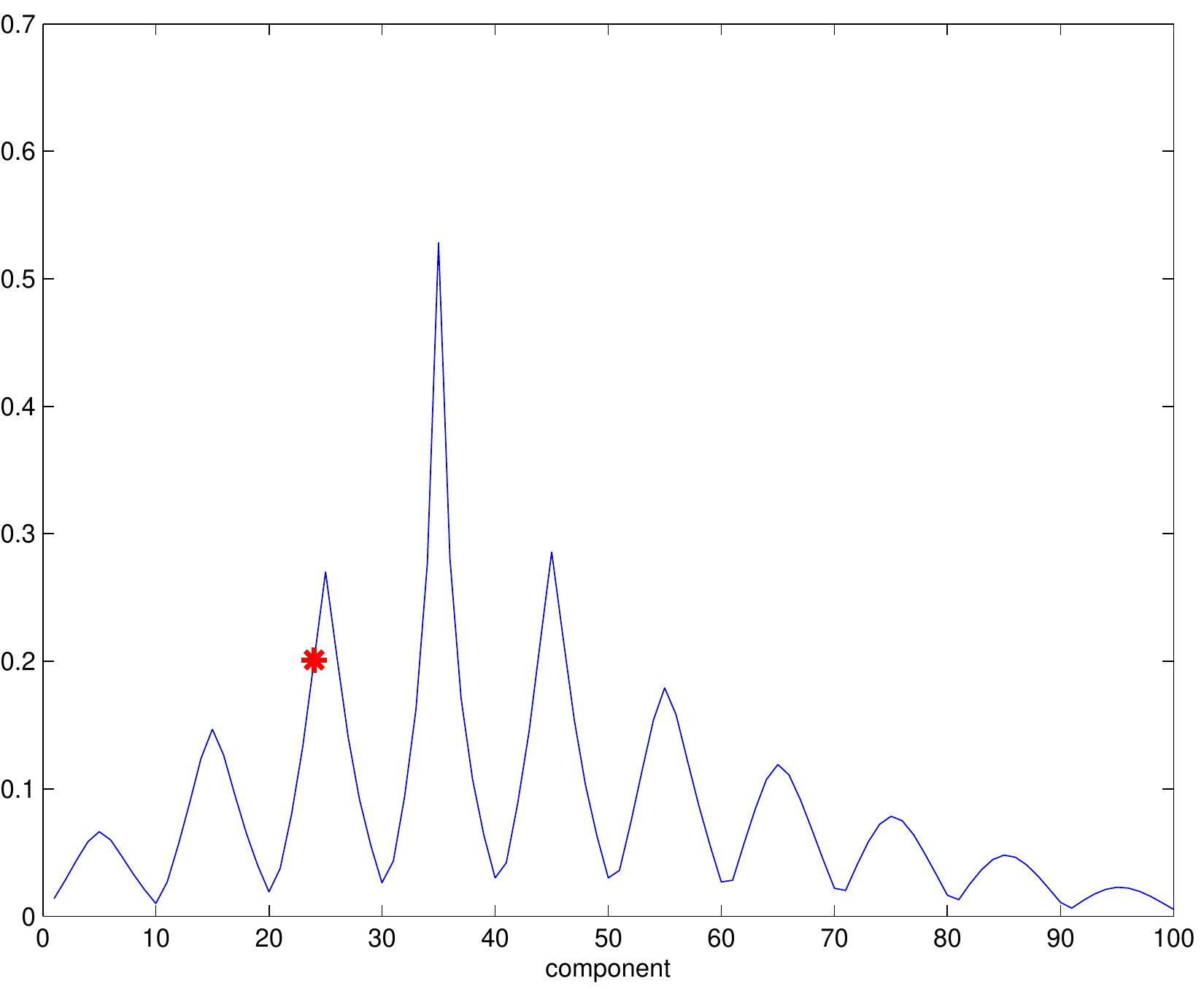} \hskip 0.5in
\includegraphics[width=1.35in,height=1.35in]{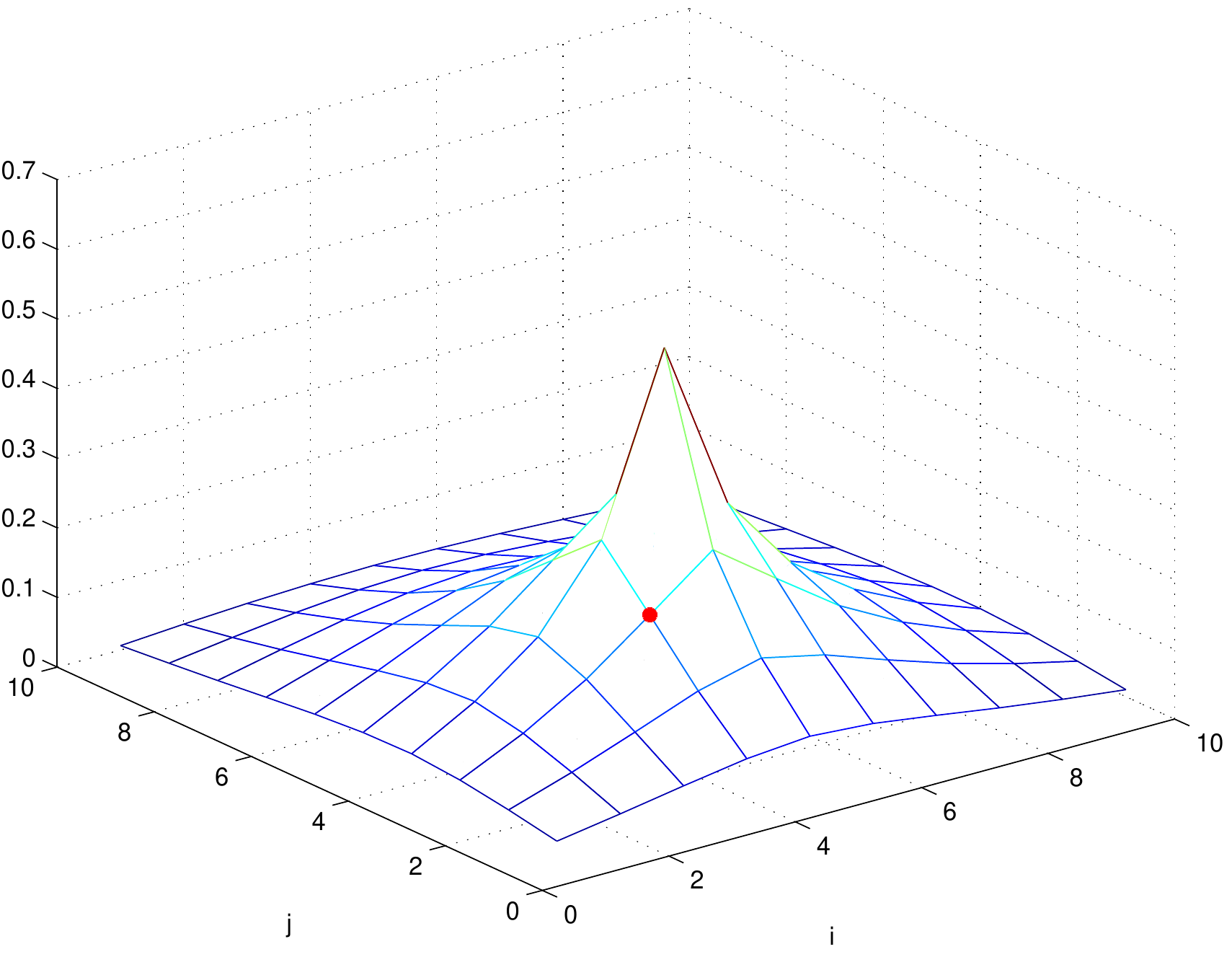} \\
\includegraphics[width=1.35in,height=1.35in]{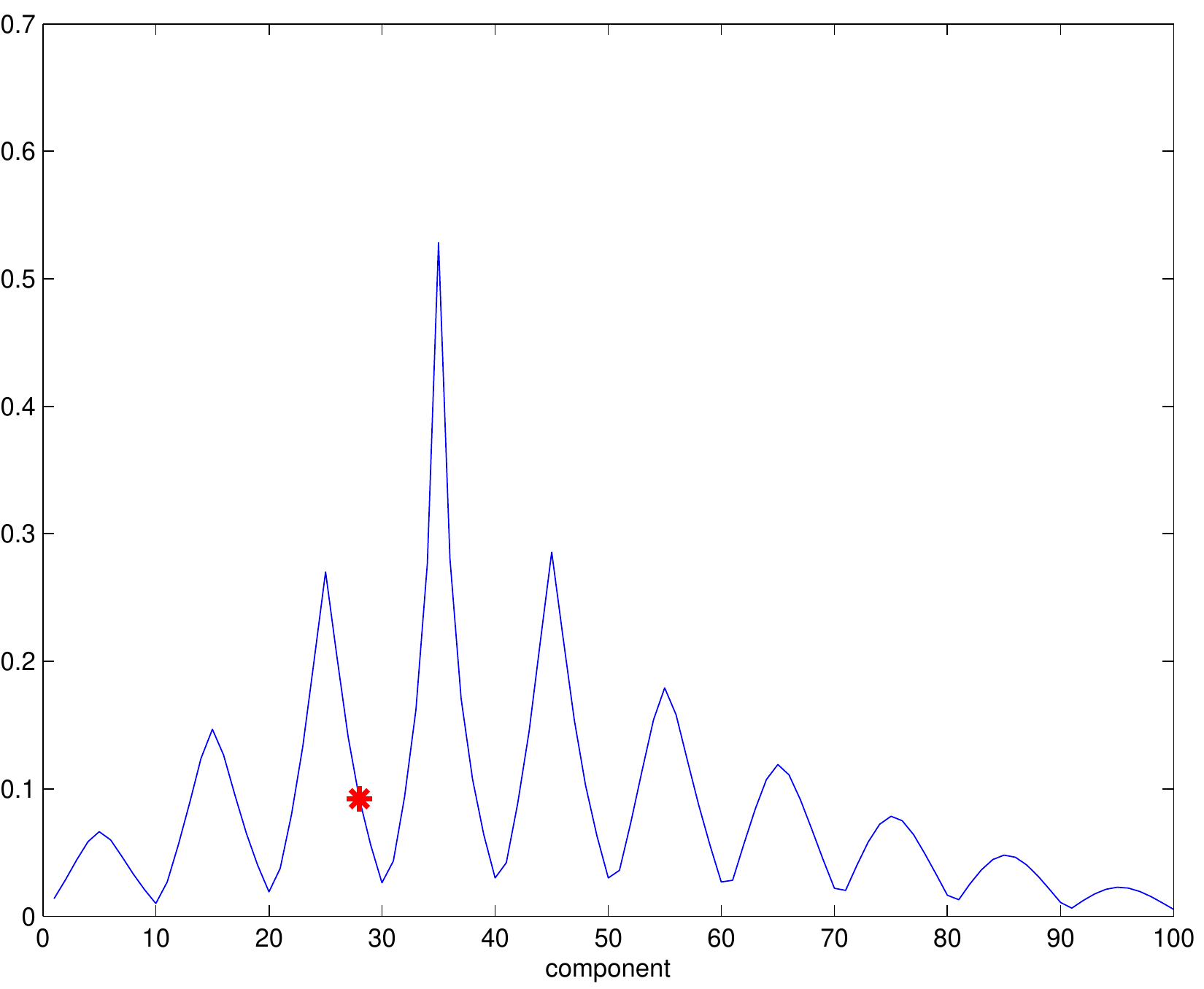} \hskip 0.5in
\includegraphics[width=1.35in,height=1.35in]{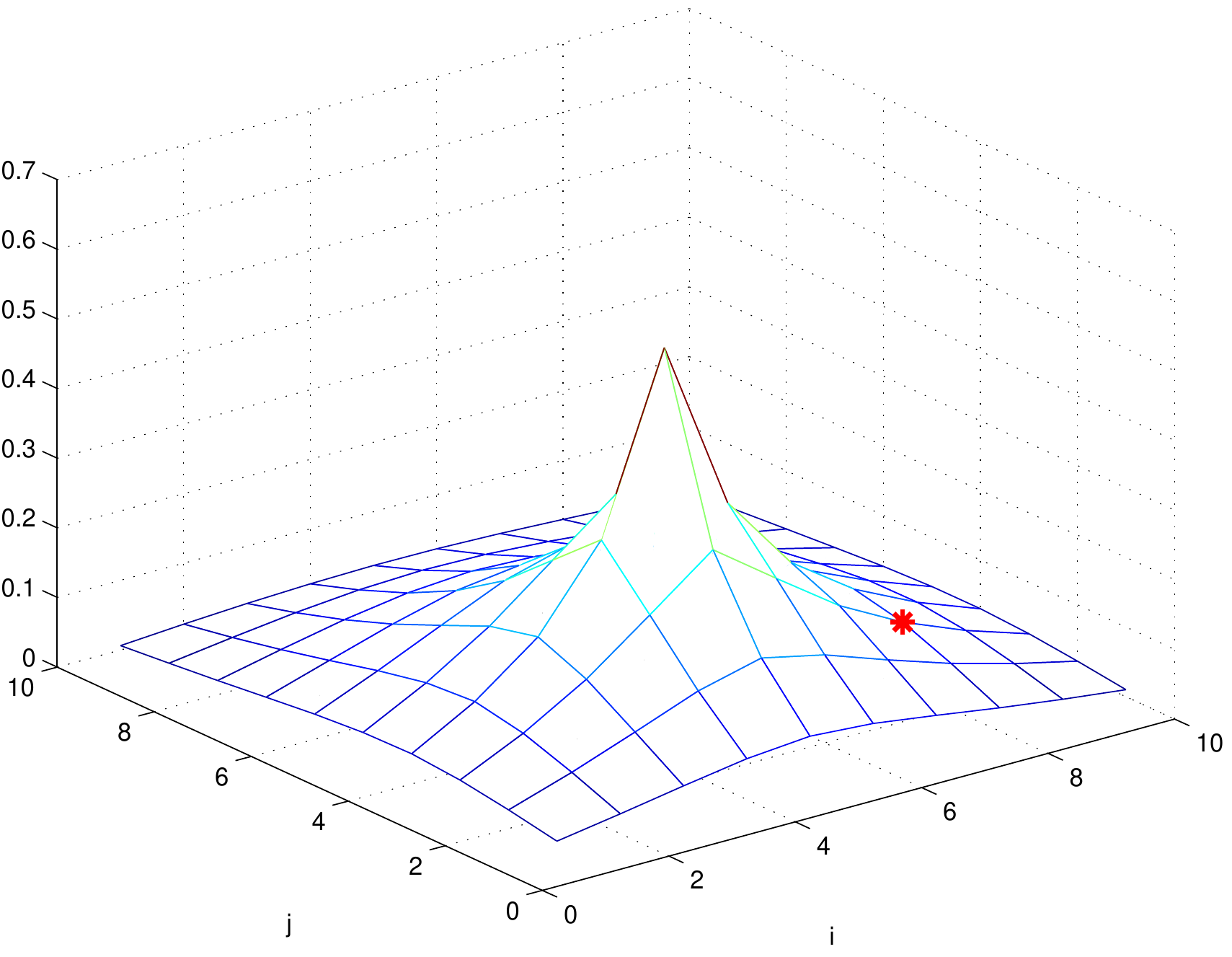} \\
\includegraphics[width=1.35in,height=1.35in]{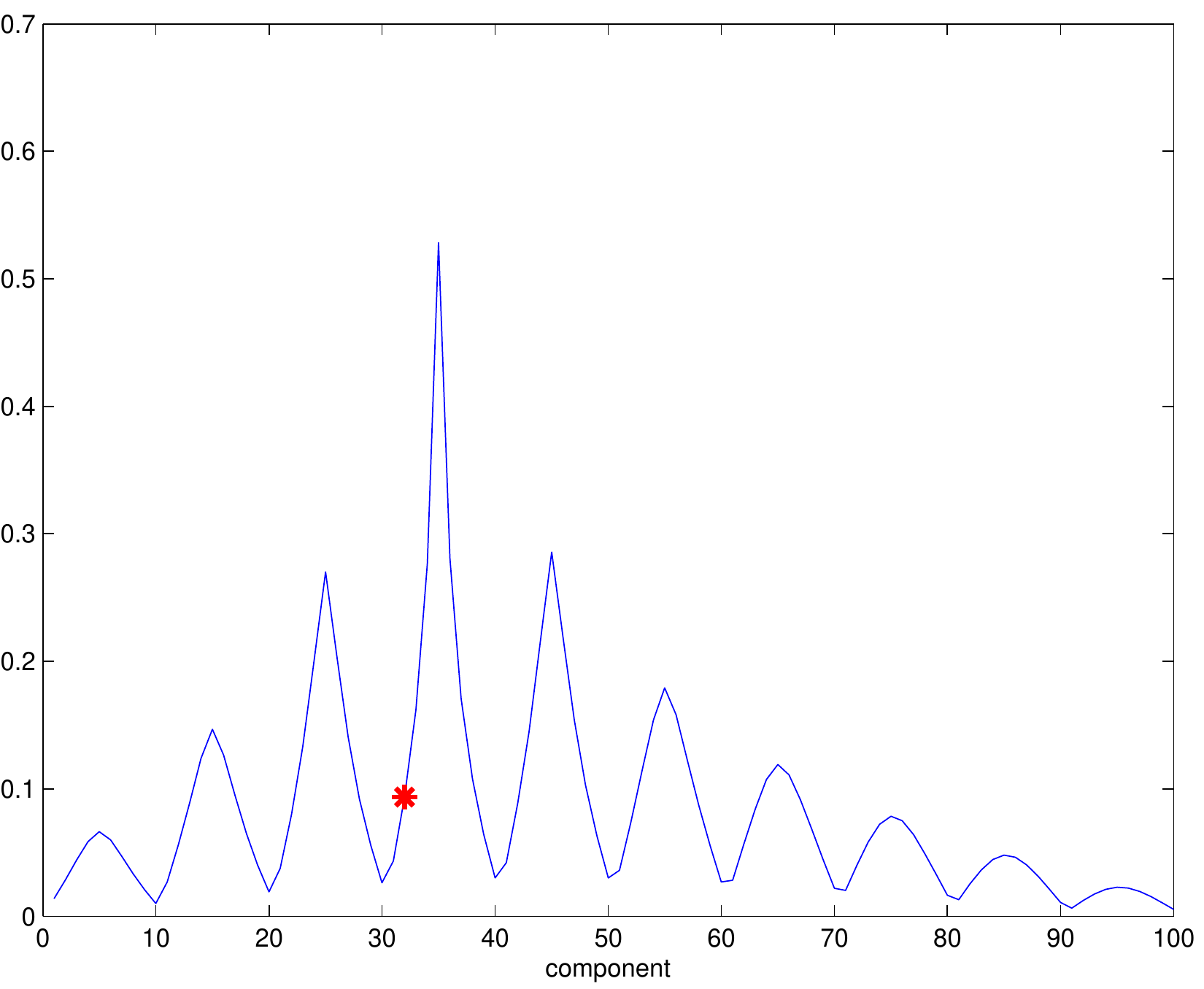} \hskip 0.5in
\includegraphics[width=1.35in,height=1.35in]{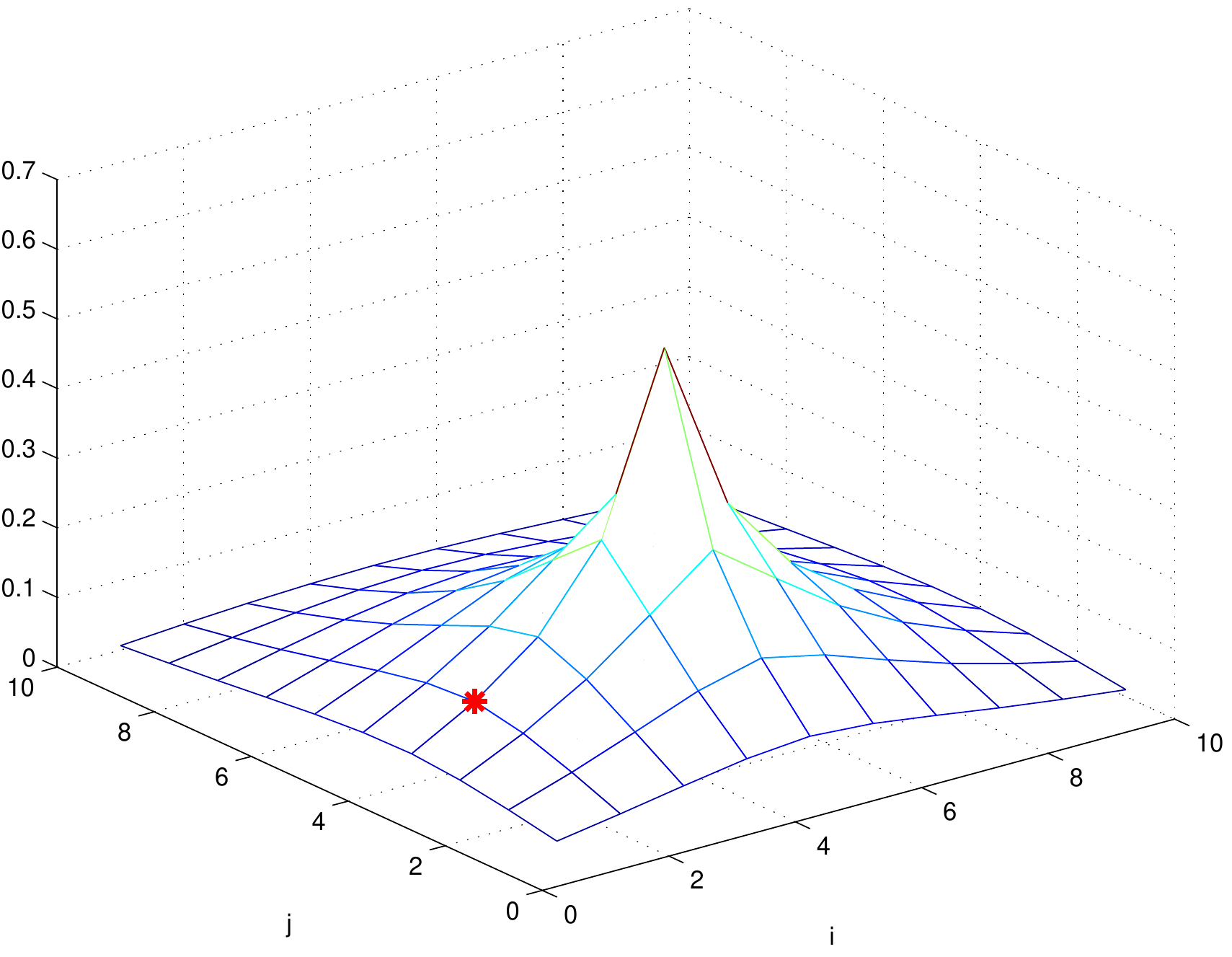} \\
\includegraphics[width=1.35in,height=1.35in]{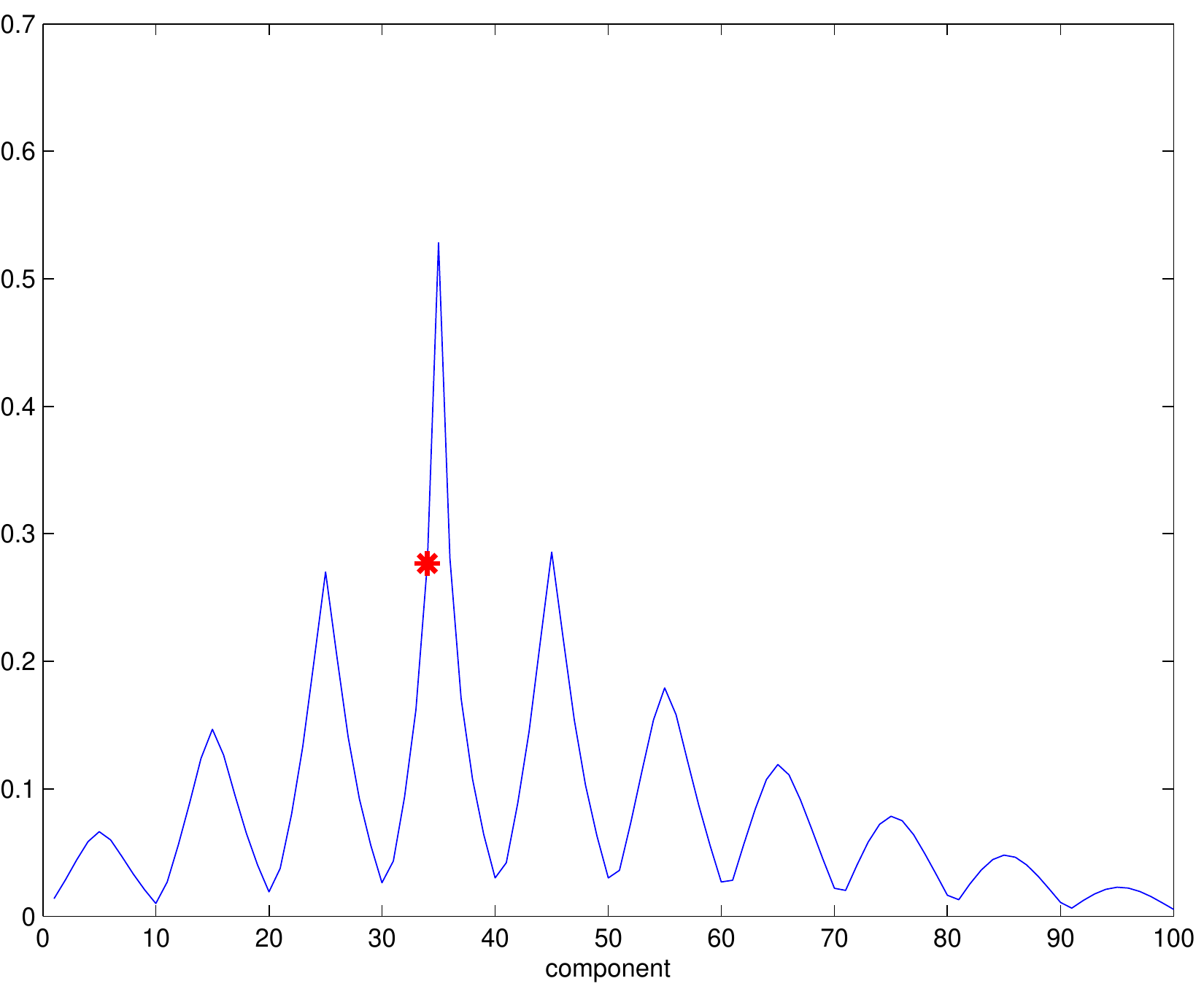} \hskip 0.5in
\includegraphics[width=1.35in,height=1.35in]{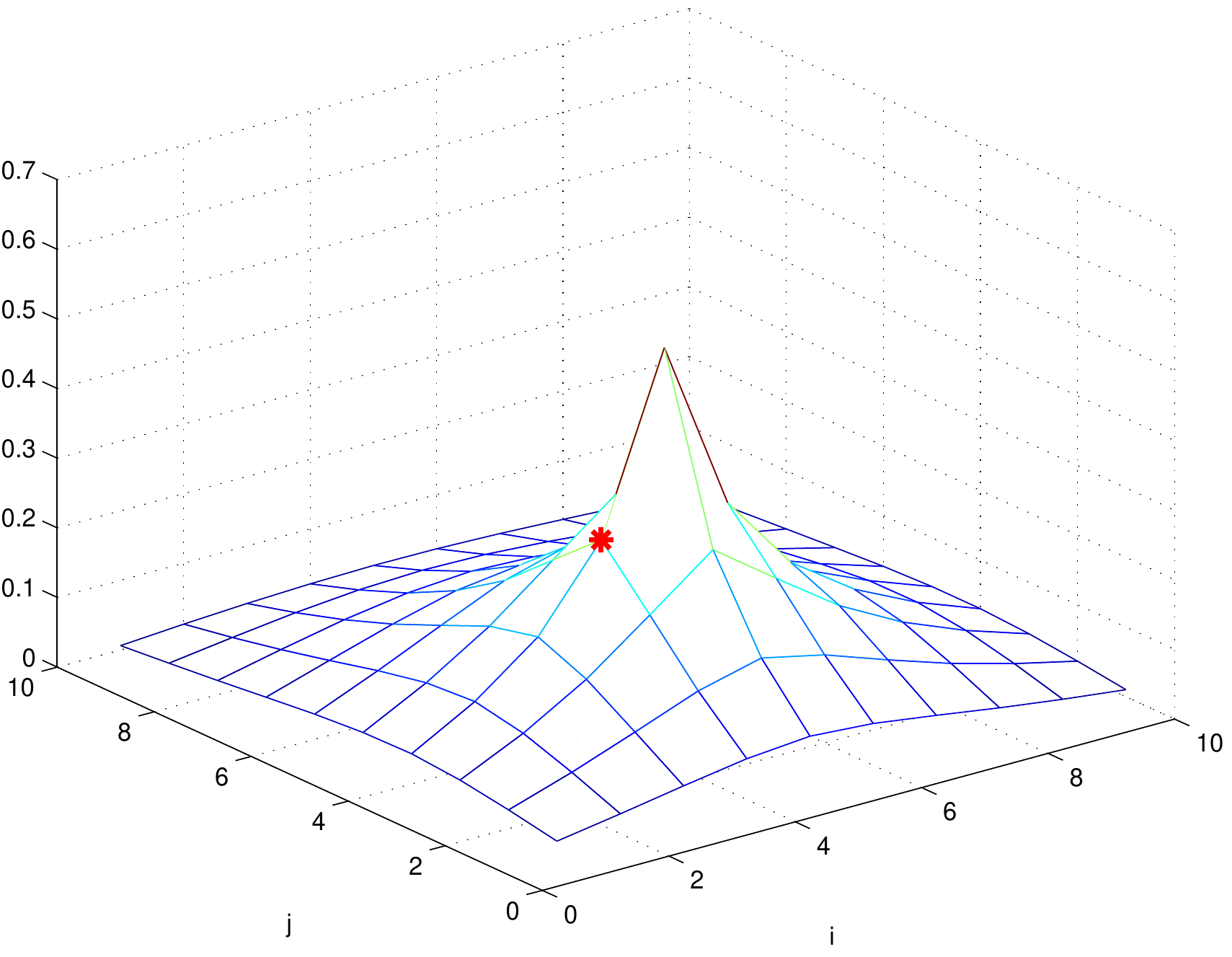} 
\caption{Column of $X$ and corresponding element on the grid. \label{fig:movie}}
\end{figure}

\begin{example}
{\rm
Let $A = {\rm tridiag}(-1,2,-1)\in\RR^{n\times n}$, $n=10$ and let $t=35$, so that
$D=e_{t_1}e_{t_2}^\TT$, $t_1=5, t_2=4$.
Figure~\ref{fig:movie} shows the connection between the column
${\cal A}^{-1} e_t$ (left plot)  and the solution $X$ (right plot on the mesh)
 of the corresponding problem
$A X + X A = e_{t_1}e_{t_2}^\TT$. The plots illustrate the correspondence of each
component of the vector as an entry of the matrix $X$, thus describing 
the oscillating behavior of the column ${\cal A}^{-1} e_t$. Note that
since $\cal A$ is also banded with bandwidth $\beta=10$, 
the overall decay is justified by the estimate in (\ref{eqn:Demko}), while
the Kronecker structure is responsible for the oscillations.
}
\end{example}

\subsection{Decay pattern of the solution matrix}
Assume the right-hand side consists of a very sparse matrix,
as an extreme case, having a single nonzero element, so that
$$
A X +  X A = {\cal E}_t, \quad {\cal E}_t = e_i \gamma e_j^\TT , \quad \gamma \ne 0.
$$
Without loss of generality, in the following we shall assume that $\gamma=1$.
Then using the integral closed form in (\ref{eqn:integral}) we can write
\begin{eqnarray}\label{eqn:X}
X = 
\frac{1}{2\pi}
\int_{-\infty}^{\infty} (\imath \omega I + A)^{-1} e_i e_j^\TT
(\imath \omega I + A)^{-\HH} {\rm d}\omega
=
\frac{1}{2\pi}
\int_{-\infty}^{\infty} z_i z_j^\HH {\rm d}\omega , 
\end{eqnarray}
where $z_i= (\imath \omega I + A)^{-1} e_i$.
Let $\mathbb S$ be the class of $n\times n$ $\beta$-banded matrices
$S$ of the form $S= \alpha_1 I + \alpha_2 S_0$
with $S_0=S_0^\HH$ (Hermitian) and $\alpha_1, \alpha_2 \in\CC$, and denote by
$[\lambda_{1}, \lambda_{2}]$ the line segment containing the eigenvalues of 
a matrix $S \in{\mathbb S}$.
Note that if $A\in{\mathbb S}$, also the matrix $\omega \imath I +A$ in
(\ref{eqn:X}) belongs to $\mathbb S$. The set ${\mathbb S}$ includes
real symmetric matrices (for $\alpha_1=0$ and $\alpha=1$), but also complex skew-Hermitian
matrices ($\alpha=0$ and $\alpha=i$), and complex shifted matrices, as they arise,
for instance, in the discretization of the Helmoltz equation. All matrices in
${\mathbb S}$ are normal matrices, that is for $S\in{\mathbb S}$ it holds
that $SS^\HH=S^\HH S$.

To be able to characterize the decay pattern in the solution $X$, we recall a
result derived in \cite{Freund1989a}.

\begin{theorem}\label{th:freund}
Let $A\in{\mathbb S}$, 
$a=(\lambda_{2}+\lambda_{1})/ (\lambda_{2}-\lambda_{1})$ and
$R>1$ be defined as
$R=\alpha + \sqrt{\alpha^2-1}$,
with $\alpha=(|\lambda_1|+|\lambda_2|)/|\lambda_2-\lambda_1|$.
Then
$$
|e_\ell^\top A^{-1}e_i| \le \frac{2R}{|\lambda_1-\lambda_2|}
B(a) \left (\frac{1}{R}\right )^{\frac{|\ell -i|}{\beta}}, \quad \ell\ne i,
$$
where, writing
$a=\alpha_R\cos(\psi)+\imath \beta_R\sin(\psi)$,
$$
B(a) := \frac {R}{   \beta_R\sqrt{\alpha_R^2-\cos^2(\psi)}(\alpha_R +
\sqrt{\alpha_R^2-\cos^2(\psi)})},
$$
with $\alpha_R=\frac 1 2 (R+\frac 1 R)$ and
$\beta_R=\frac 1 2 (R-\frac 1 R)$.
\end{theorem}

With this result, the following bound for the entries of the
Lyapunov solution matrix can be obtained; in spite of the more
general setting, the proof proceeds as in
\cite{CanutoSimonciniVeraniLAA.14}, and it is therefore omitted.

\begin{theorem}\label{th:entries}
For $\omega\in\RR$, 
let $\omega\imath I + A \in{\mathbb S}$, with eigenvalues contained in the
line segment $[\lambda_1, \lambda_2 ] := 
[\lambda_{\min}(A)+\omega\imath, \lambda_{\max}(A)+\omega\imath]$.
With the notation of Theorem \ref{th:freund},
and
 for $k=(k_1,k_2)$ and $t=(t_1, t_2)$, $k_i,t_i\in\{1, \ldots, n\}$, $i=1,2$,
%
the following holds.

i) If $t_1\ne k_1$ and $t_2\ne k_2$, then
$$
|(X)_{k_1,k_2}|  
\le
\frac{1}{2\pi}
\frac{64}{|\lambda_{2}-\lambda_{1}|^2} \int_{-\infty}^{\infty} 
 \left ( \frac{R^2}{(R^2-1)^2}\right )^2
\left ( \frac{1}{R}\right )^{|t_1-k_1|/\beta+|t_2-k_2|/\beta-2} {\rm d}\omega ;
$$
ii) If either $t_1= k_1$  or $t_2= k_2$, then
$$
|(X)_{k_1,k_2}|  
\le
\frac{1}{2\pi}
\frac{8}{|\lambda_{2}-\lambda_{1}|} \int_{-\infty}^{\infty} 
\frac 1 {\sqrt{|\lambda_{1}|^2 + \omega^2}}
 \frac{R^2}{(R^2-1)^2}
\left ( \frac{1}{R}\right )^{|t_1-k_1|/\beta +|t_2-k_2|/\beta-1} {\rm d}\omega ;
$$
iii) If both $i= \ell$  and $j= m$, then
$$
|(X)_{k_1,k_2}|  
\le
\frac{1}{2\pi}
\int_{-\infty}^{\infty} \frac 1 {|\lambda_{1}|^2 + \omega^2} {\rm d}\omega = \frac 1 {2|\lambda_{1}|}.
$$
\end{theorem}

\begin{figure}[htb]
\centering
\includegraphics[width=1.5in,height=1.5in]{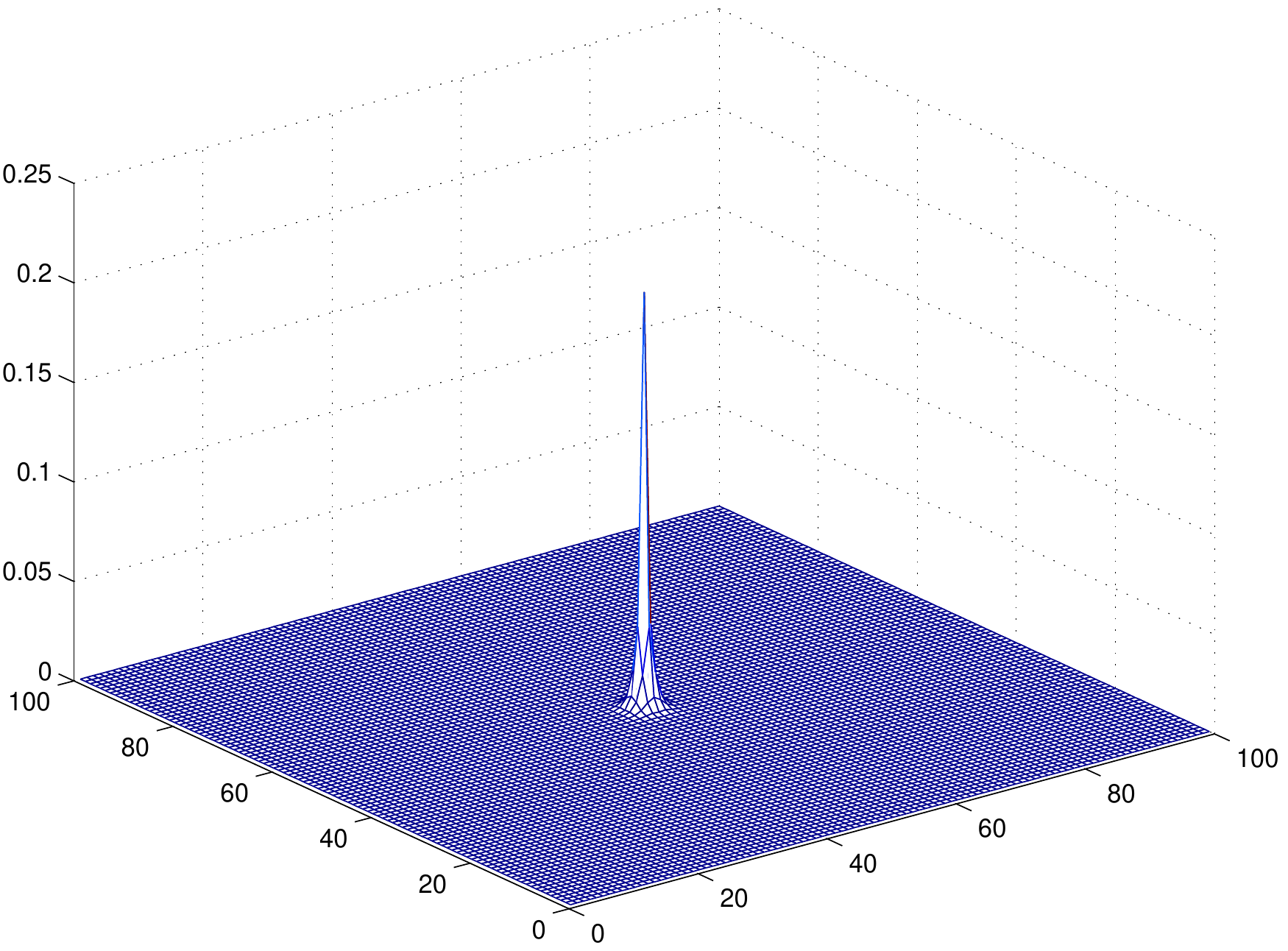} 
\hskip 0.5in
\includegraphics[width=1.5in,height=1.5in]{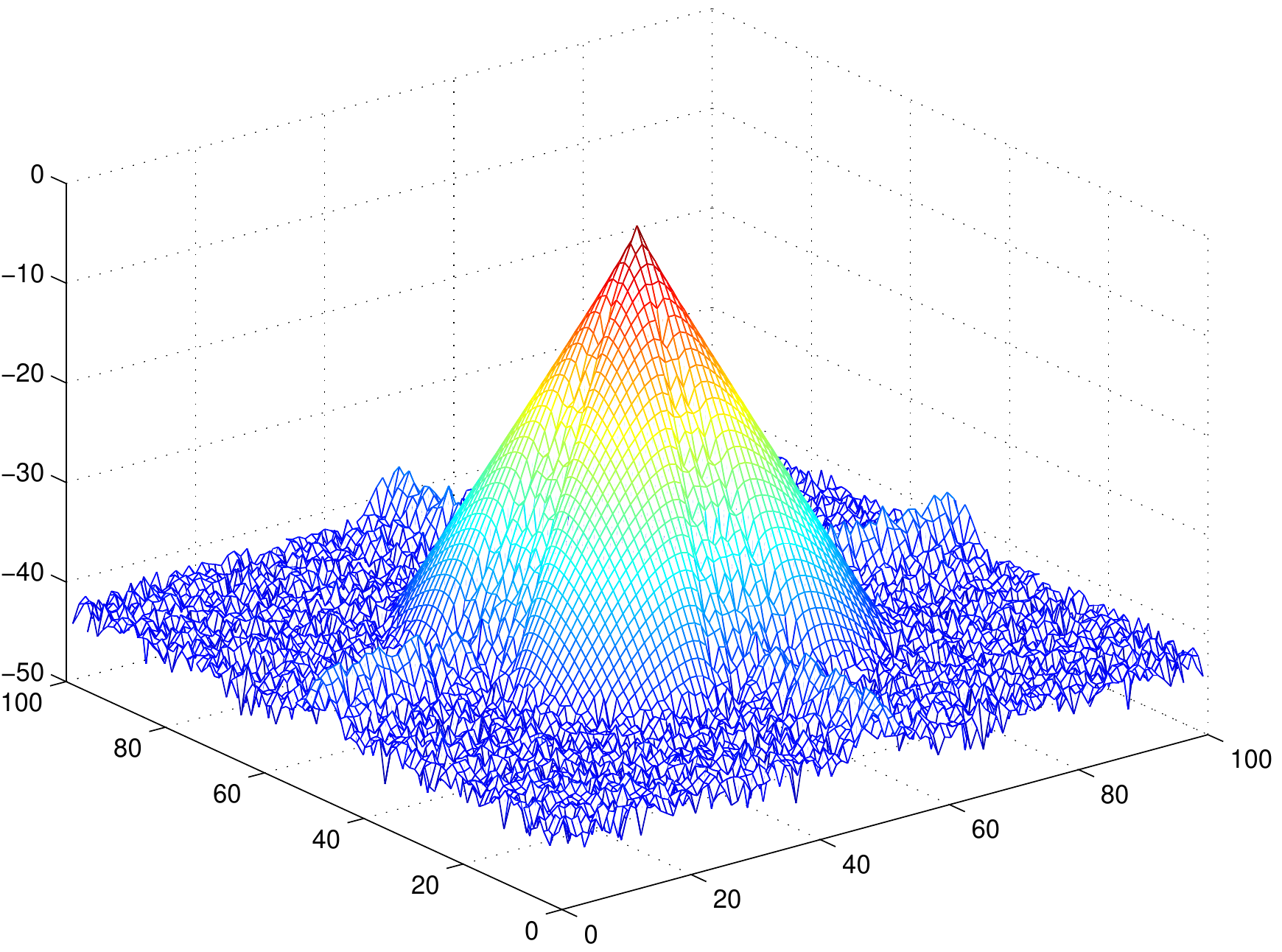}
\caption{Decay pattern of $X$ for Example \ref{ex:Skew}.
Left: linear scale. Right: logarithmic scale. \label{fig:Skew}}
\end{figure}

\begin{example}\label{ex:Skew}
{\rm
We consider the $n\times n$ nonsymmetric tridiagonal matrix 
$M={\rm tridiag}(1,2,-1)$, $n=100$, 
 and $D=bb^\TT$ with $b=e_{50}$; note that the same pattern is obtained
with its complex counterpart $M={\rm tridiag}(i,2,i) \in {\mathbb S}$.
Figure~\ref{fig:Skew} shows the solution
pattern in linear (left) and logarithmic (right) scale.
The solution entries decay very quickly away from the entry $(50,50)$. This fact
can be appreciated by noticing the $z$-axis in the logarithmic scale.
}
\end{example}
\vskip 0.1in

In cases of interest in applications, such as in the case of
problems stemming from certain discretizations of elliptic partial differential
equations, the matrix
$A$ itself may have  a structure of sum of Kronecker products, that is
 \begin{eqnarray}\label{eqn:Akron}
A=M\otimes I + I \otimes M,
\end{eqnarray}
where $M$ is $\sqrt{n}\times \sqrt{n}$ and banded.
We use once again (\ref{eqn:X}) for expressing $e_{k_1}^\TT X e_{k_2}$, by first
assuming that a lexicographic order was used to generate the indexes $k_1$, $t_1$ on
the two-dimensional grid, so that we can identify
$k_1=(k_{11},k_{12})$, $t_1=(t_{11},t_{12})$.
Then we can see that the inner product inside the integral satisfies
$e_{k_1}^\TT(\imath \omega I + A)^{-1} e_{t_1} = e_{k_{11}}^\TT Z_{t_1} e_{k_{12}}$
where $Z_{t_1}$ solves the Lyapunov equation
$$
\left (M+\frac{1}{2}\omega \imath I\right) Z + 
Z \left (M+\frac{1}{2}\omega\imath I\right ) = e_{t_{11}} e_{t_{12}}^\TT.
$$
Therefore, in a recursive manner, the decay of the entries in the matrix $Z$ can 
be described in terms of Theorem \ref{th:entries}.
In practice, this pattern results in a ``local'' oscillation associated with the
decay of the columns of $Z$, and a ``global'' oscillation, associated with the
decay of $X$. 

\begin{figure}[htb]
\centering
\includegraphics[width=1.4in,height=1.4in]{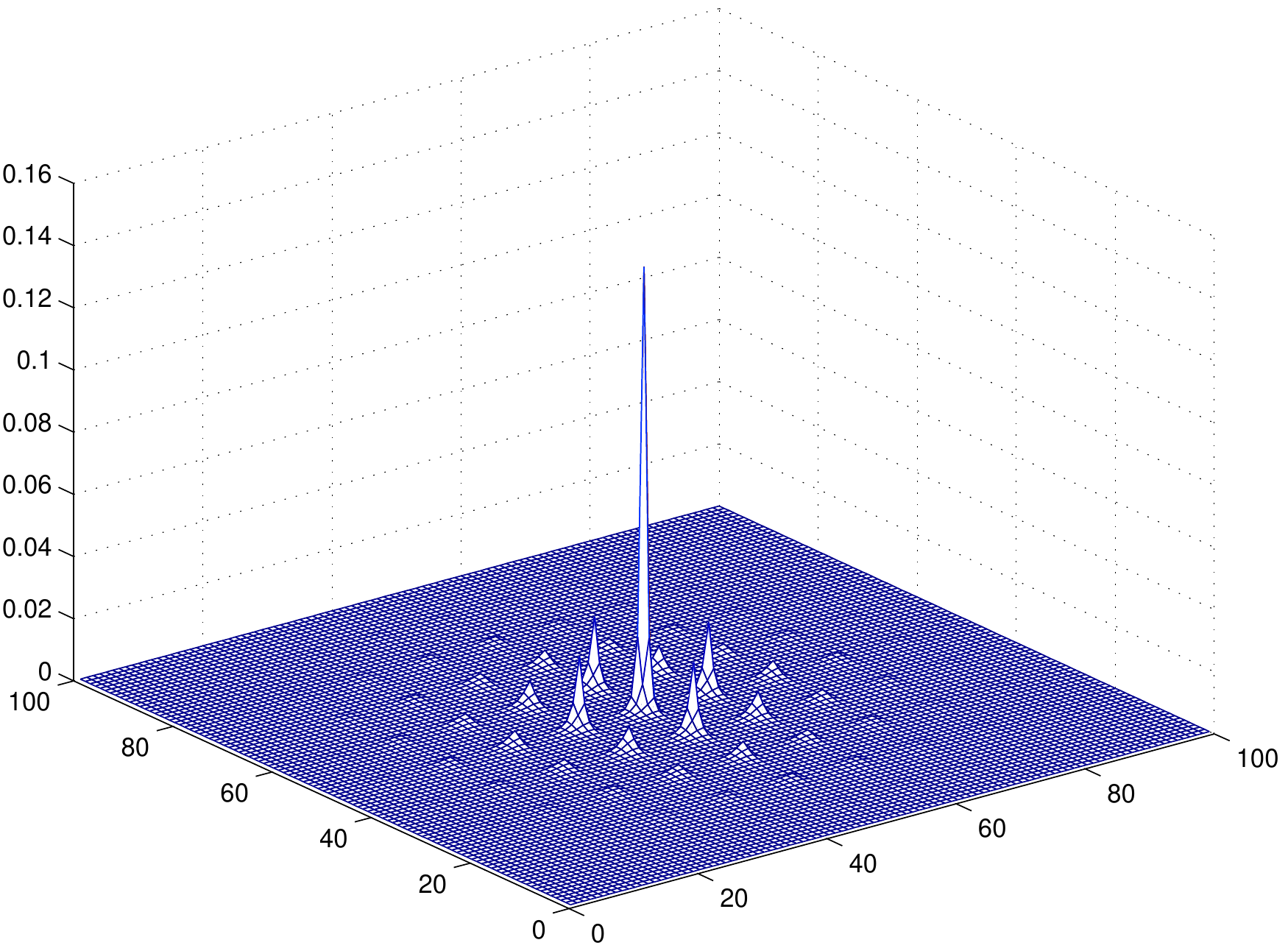}\hskip 0.2in
\includegraphics[width=1.4in,height=1.4in]{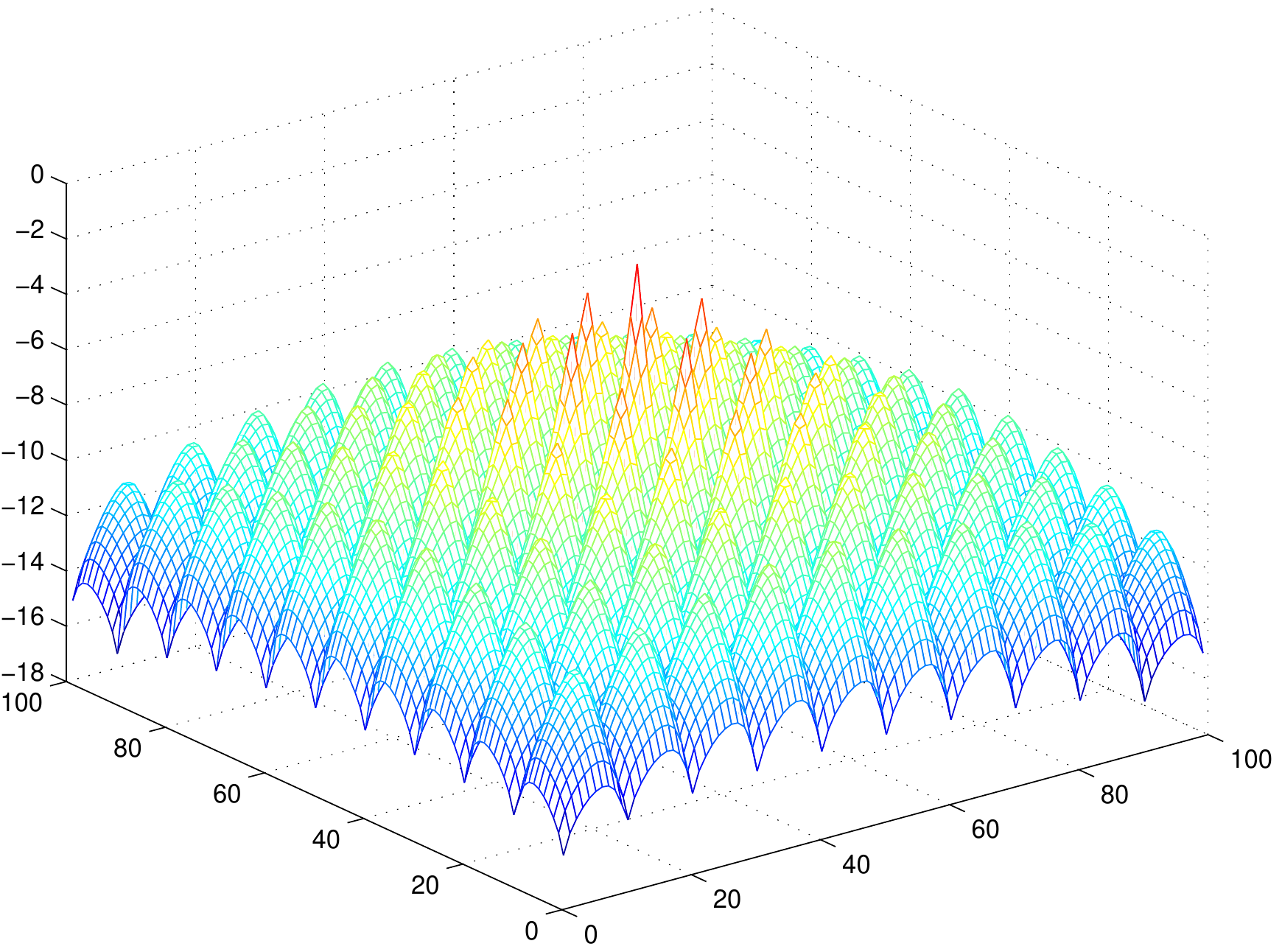}\hskip 0.2in
\includegraphics[width=1.4in,height=1.4in]{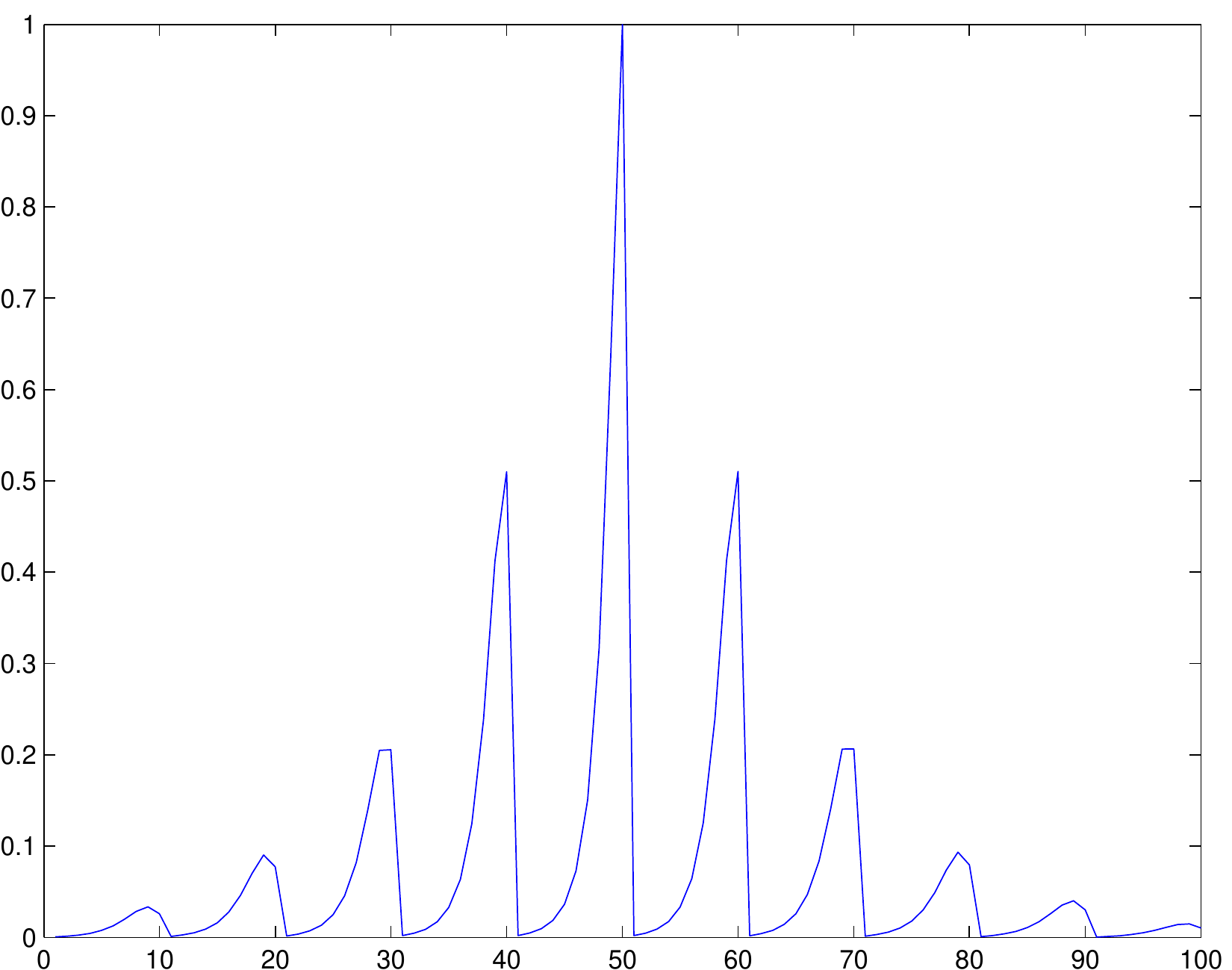}
\caption{Pattern of solution $X$ with $A$ in (\ref{eqn:Akron}). Left: linear scale. 
Middle: logarithmic scale. Right: $48$th column of $X$ (slice of left plot in linear scale).
\label{fig:3D}}
\end{figure}

\vskip 0.1in
\begin{example}
{\rm
A typical sample of this behavior 
is shown in Figure \ref{fig:3D}, where 
$M={\rm tridiag}(-1,2,-1)\in\RR^{\sqrt{n}\times \sqrt{n}}$, $\sqrt{n}=10$, 
so that $A\in\RR^{n\times n}$, and $B=e_{50}$.
}
\end{example}
\vskip 0.1in

Combining the Kronecker structure of $A$ with the one in (\ref{eqn:Akron})
 of the Lyapunov equation, we can reformulate the problem as
$$
(M\otimes I \otimes I \otimes I + I\otimes M \otimes I \otimes I +
I\otimes I \otimes M \otimes I + I\otimes I \otimes I \otimes M) x = b,
$$
where $b= {\rm vec}(BB^\TT)$, and
all identity matrices have dimension $\sqrt{n}$. This form reveals the actual $4\times 4$
tensorial nature of the solution $X$ in case $A$ has the described
 Kronecker structure. Computational methods to solve for $x$ while using only
information in $\RR^{\sqrt{n}}$ together with the tensor structure
 have been proposed; see, e.g., \cite{Kressner.Tobler.10}.
Here we just want to emphasize that the quality of the decay in the entries of $X$, and thus of $x$,
is all determined by that of the banded matrix $M$, while the location of
this decay pattern is ruled by the Kronecker form of $A$, as nicely
visible in the middle plot of Figure~\ref{fig:3D}.

\subsection{Numerical low rank properties of the solution matrix}
In the literature there is very large experimental evidence that
if the right-hand side in (\ref{eqn:Lyap}) has low rank, then the eigenvalues
of the solution $X$  decay very rapidly. In some cases, such
as for $A$ symmetric and positive definite, this behavior has been
proved \cite{Penzl2000}. Bounds that cope with non-normality have also
been derived, see, e.g., \cite[sec.3.1.2]{Sabino.PhD.06} and \cite{ASZ.02}.
Theoretical results showing the fast decay of the spectrum rely
on the low rank of the matrix $D$, and on the spectral properties of
the coefficient matrix $A$. We are not aware of results that combine
the truncation strategies associated with
the matrix decay pattern with the low rank property. Both strategies aim
at drastically reducing the memory requirements of the solution matrix, and in general,
sparsity-driven truncated matrices may have larger numerical rank than the original 
matrices. However, the following example provides an interesting setting,
which will be worth exploring in the future.

\begin{figure}[htb]
\centering
\includegraphics[width=1.6in,height=1.6in]{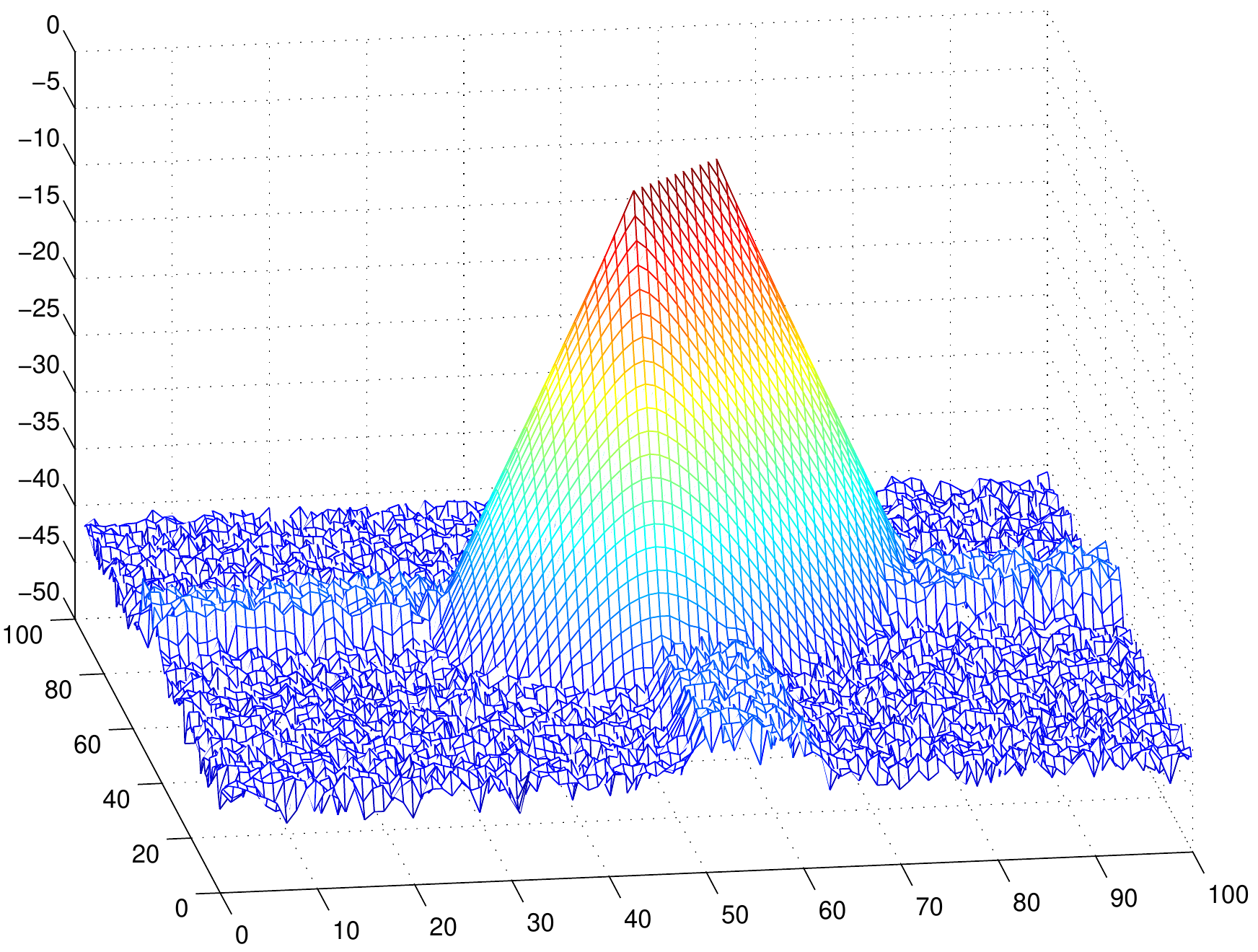}\hskip 0.2in
\includegraphics[width=1.6in,height=1.6in]{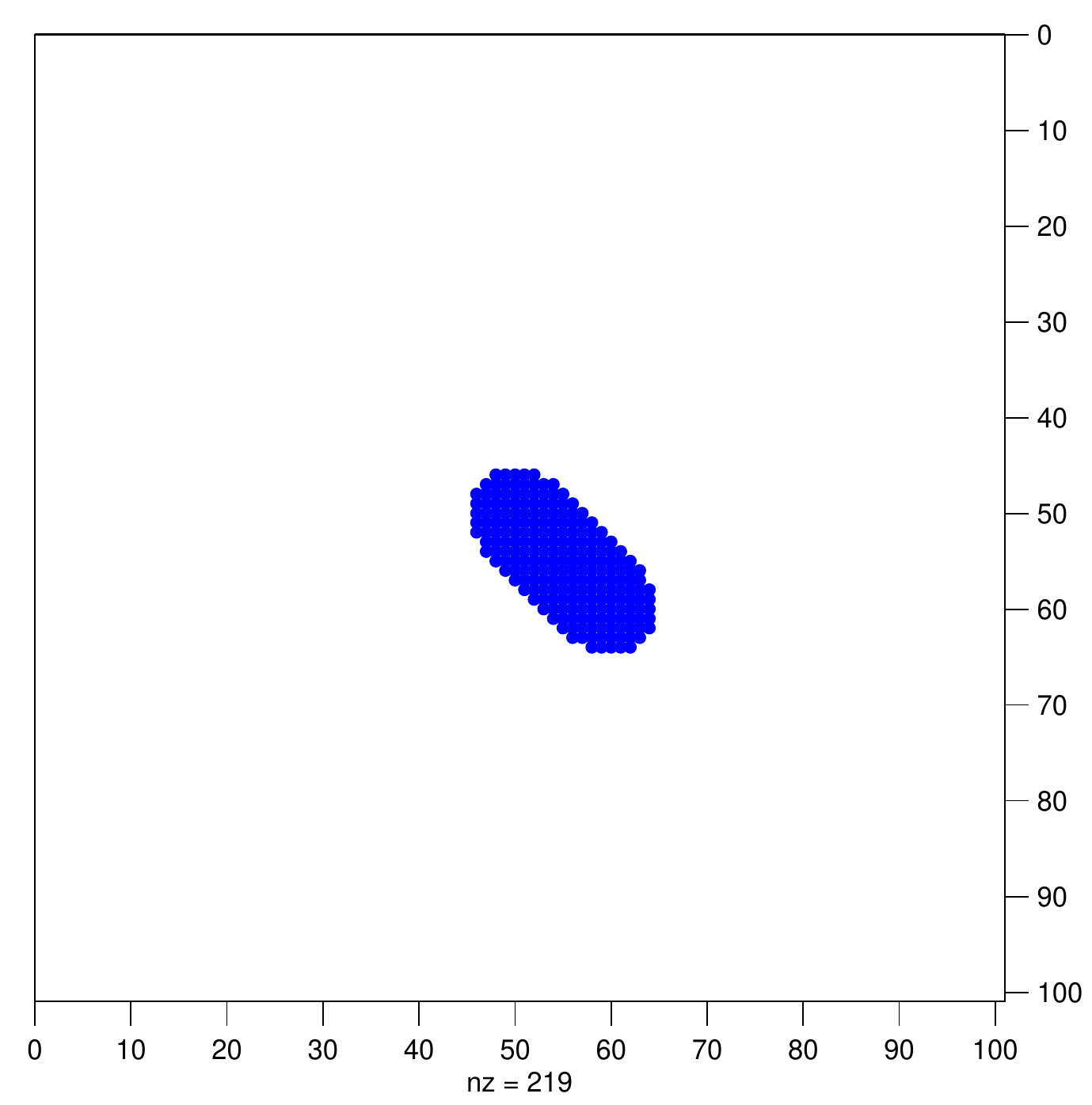}
\caption{Solution $X$ for Example \ref{ex:trunc}. Left: pattern of $X$ with
logarithmic scale, {\tt nnz}($X$) = 9724. Right: Sparsity pattern of
truncated version of $X$: all entries below $10^{-5}$ are omitted.
\label{fig:trunc}}
\end{figure}

\begin{example}\label{ex:trunc}
{\rm
Let $A={\rm tridiag}(-1,4,-1)\in\RR^{n\times n}$, $n=100$
 and $D=BB^\TT$ with $B=[e_{50},\ldots, e_{60}]$.
The left plot of Figure~\ref{fig:trunc}
shows the solution pattern. The symmetric and positive
semidefinite matrix $X$ has 9724 nonzero elements. Moreover, 
25 eigenvalues of $X$ are above $10^{-14}$, therefore
$X$ can be very precisely reproduced as $X=X_1X_1^\TT$ with a
tall matrix $X_1$ having rank (and thus number of columns) as low as 25.
The right plot shows the sparsity pattern of the matrix $\widetilde X$
obtained by zeroing all entries of $X$ below $10^{-5}$ (the largest entry
of $X$ is ${\cal O}(1)$). The matrix $\widetilde X$ has precisely 19
nonzero singular values, and only  219 nonzero entries (see right plot
of Figure~\ref{fig:trunc}). The accuracy of
$\widetilde X$ is not as good as with the rank-25 matrix, since
$\|X-\widetilde X\| \approx 10^{-5}$, however for most applications
this approximation will be rather satisfactory.
Though this example is not sufficiently general to make a general
case, this phenomenon deserves further analysis.
}
\end{example}

\section{Numerical considerations}
Sparsity and quasi-sparsity properties are now starting to be considered 
in the numerical approximation of the solution to (\ref{eqn:Lyap}).
In the following we shall discuss two situations where the
analysis of the previous sections can provide insights into
the understanding and development of numerical strategies.

\subsection{Decay issues in the approximate solution by projection}
Let $A$ be $n\times n$ and sparse, $D=BB^\TT$ with $B\in\RR^{n\times s}$, $s\ll n$.
An effective strategy for approximately solving 
(\ref{eqn:Lyap}) consists of determining a good approximation 
space $K_m$ of small dimension $m$, and then generating $X_m=V_m Y V_m^\TT 
\approx X$ where the orthonormal columns of $V_m\in\RR^{n\times m}$ span $K_m$.
The matrix $Y\in\RR^{m\times m}$ is determined by imposing some
extra condition. For instance, a possible strategy requires that
the residual $R_m := A X_m + X_m A^\TT - BB^\TT$ satisfies the
following (orthogonality) condition:
$$
V_m^\TT R_m V_m = 0.
$$
Inserting the expression for the residual, and recalling that
$V_m^\TT V_m = I_m$, this condition leads
to the following reduced Lyapunov equation 
\begin{eqnarray}\label{eqn:reduced}
(V_m^\TT A V_m) Y + Y (V_m^\TT A^\TT V_m) = V_m^\TT B B^\TT V_m .
\end{eqnarray}
Setting $T_m = V_m^\TT A V_m$ and $\widehat B_m = V_m^\TT B$, the
matrix $Y$ is thus determined by solving the equation
$T_m Y + Y T_m^\TT = \widehat B_m \widehat B_m^\TT$ of (small) size $m$;
this can be performed by using a ``dense'' method, such as the popular
Bartels-Stewart algorithm, whose computational cost is ${\cal O}(m^3)$ \cite{Bartels.Stewart.72}.
The effectiveness of the whole approximation thus depends on how
small $m$ can be and still obtain a satisfactory approximation.

The choice of $K_m$ is crucial to obtain a rich enough approximation
space while keeping its dimension $m$ small.
We refer the reader to the recent survey \cite{Simoncini.survey13} for
a comprehensive description of these issues and the relevant recent
literature. Here we consider the
simplest possible choice of approximation space in the class of
Krylov subspace methods: assuming for ease of presentation
$B=b\in\RR^n$, we define the (standard) Krylov subspace as the
vector space
$$
K_m =K_m(A,b) :=  {\rm span}\{b, Ab, \ldots, A^{m-1}b\} .
$$
As $m$ increases, $K_{m+1}$ is obtained from $K_m$ by one additional
multiplication by $A$ of the last basis vector. An orthonormal basis
can be progressively constructed as $m$ increases
by using the Arnoldi procedure \cite{Saad2003}, giving rise to the
columns of $V_m$; in particular, $b=V_m e_1 \|b\|$. The Arnoldi recurrence
also establishes the following relation
\begin{eqnarray}\label{eqn:Arnoldi}
A V_m = V_m T_m + v_{m+1} t_{m+1} e_m^\TT ,
\end{eqnarray}
where $v_{m+1}$ is the next vector of the basis, and $t_{m+1}$ is computed
during the recurrence.

A nice feature of the Arnoldi procedure, is that the reduced
matrix $T_m = V_m^\TT A V_m$ is upper Hessenberg, so that the
cost of solving the reduced matrix equation in (\ref{eqn:reduced}) is lower
than for a full matrix. For $A$ symmetric, the Arnoldi procedure
reduces to the Lanczos procedure, and most importantly in our
context, $T_m$ is also symmetric and thus tridiagonal. With these
hypotheses, the sparsity considerations of the previous sections become
relevant. 
Indeed, the solution $Y$ to (\ref{eqn:reduced}) can be written in closed form as
$$
Y = 
 \frac 1 {2\pi} \int_{-\infty}^{\infty}
(\omega \imath I -T_m)^{-1}e_1\|b\|^2 e_1^\TT 
(\omega \imath I - T_m)^{-\HH} d\omega  ,   \quad \widehat B \widehat B^\TT = 
(e_1 \|b\|) (e_1 \|b\|)^\TT .
$$
It is also important to keep in mind that a new $Y$ will be constructed each time
the Krylov subspace is enlarged by one, since $T_m$ will be extended to $T_{m+1}$
by the addition of one row and one column.

The entries of $Y$ can be bounded as described in 
Theorem~\ref{th:entries}, where however the right-hand side this time
has a fixed nonzero entry, the one corresponding to the (1,1) position.
Thus,
\begin{eqnarray}\label{eqn:Yentries}
Y_{i,j} = Y_{j,i} =
 \frac 1 {2\pi} \int_{-\infty}^{\infty}
e_i^\TT (\omega \imath I -T_m)^{-1}e_1\|b\|^2 e_1^\TT 
(\omega \imath I - T_m)^{-\HH} e_j d\omega  ,  
\end{eqnarray}
Since $T_m$ is tridiagonal, and thus banded with bandwidth $\beta=1$, the quantities 
$|e_i^\TT (\omega \imath I -T_m)^{-1}e_1|$ undergo the exponential decay
described in Theorem \ref{th:freund}, with the rate $(1/R)^{|i-1|}$. Here $R$ is
associated with the spectral properties of $T_m$. We recall here that due to the
Courant-Fisher min-max theorem, the
eigenvalues of $T_m$ tend to approximate the eigenvalues of $A$  when $A$
is symmetric, therefore the spectral properties of $T_m$ will eventually be
associated with those of $A$, as $m$ grows.

\vskip 0.1in
\begin{remark}
Our presentation is based on (\ref{eqn:Yentries}), which exploits decay properties
of the inverses of shifted matrices described in Theorem~\ref{th:freund}.
Qualitatively similar results could be obtained by using the exponential
closed form in (\ref{eqn:cauchy}), by using very recently developed bounds
for the entries of the exponential matrix \cite{Benzi.Simoncini.15}.
\end{remark}
\vskip 0.1in

It is important to realize that the presence of a fixed column index $e_1$ in (\ref{eqn:Yentries})
provides additional structure to the decay pattern. For instance, 
all diagonal entries of $Y$ have a decreasing pattern, as
$$
Y_{i,i} =
- \frac {\|b\|^2} {2\pi} \int_{-\infty}^{\infty}
|e_i^\TT (\omega \imath I -T_m)^{-1}e_1|^2 \|b\|^2 
d\omega  ,
$$
so that, using Theorem \ref{th:entries},
$$
|(Y)_{i,i}|  
\le
\frac{\|b\|^2}{2\pi}
\frac{64}{|\lambda_{\max}-\lambda_{\min}|^2} \int_{-\infty}^{\infty} 
 \left ( \frac{R^2}{(R^2-1)^2}\right )^2
\left ( \frac{1}{R}\right )^{2|i-1|-2} {\rm d}\omega , \quad i>1 ;
$$
where $R$ is as described in Theorem \ref{th:freund}, with $T_m$ real symmetric. More
explicit estimates can be obtained by bounding the integral from above, as done, e.g.,
in \cite{CanutoSimonciniVeraniLAA.14}; we avoid these technicalities in this presentation.

The bound above allows us to distinguish between the decay of the entries of $|Y|$ and that of
the entries of $|T_m^{-1}|$, since in the latter case, the diagonal entries 
do not necessarily decay.
Moreover, the integrand for the nondiagonal entry $Y_{i,j}$, $j\ne i$  will contain the term
$$
\left ( \frac{1}{R}\right )^{|i-1|+|j-1|-2}, \qquad i,j>1 ,
$$
illustrating a superexponential decay of the antidiagonals, as $i$ grows.

\begin{figure}[htb]
\centering
\includegraphics[width=1.5in,height=1.5in]{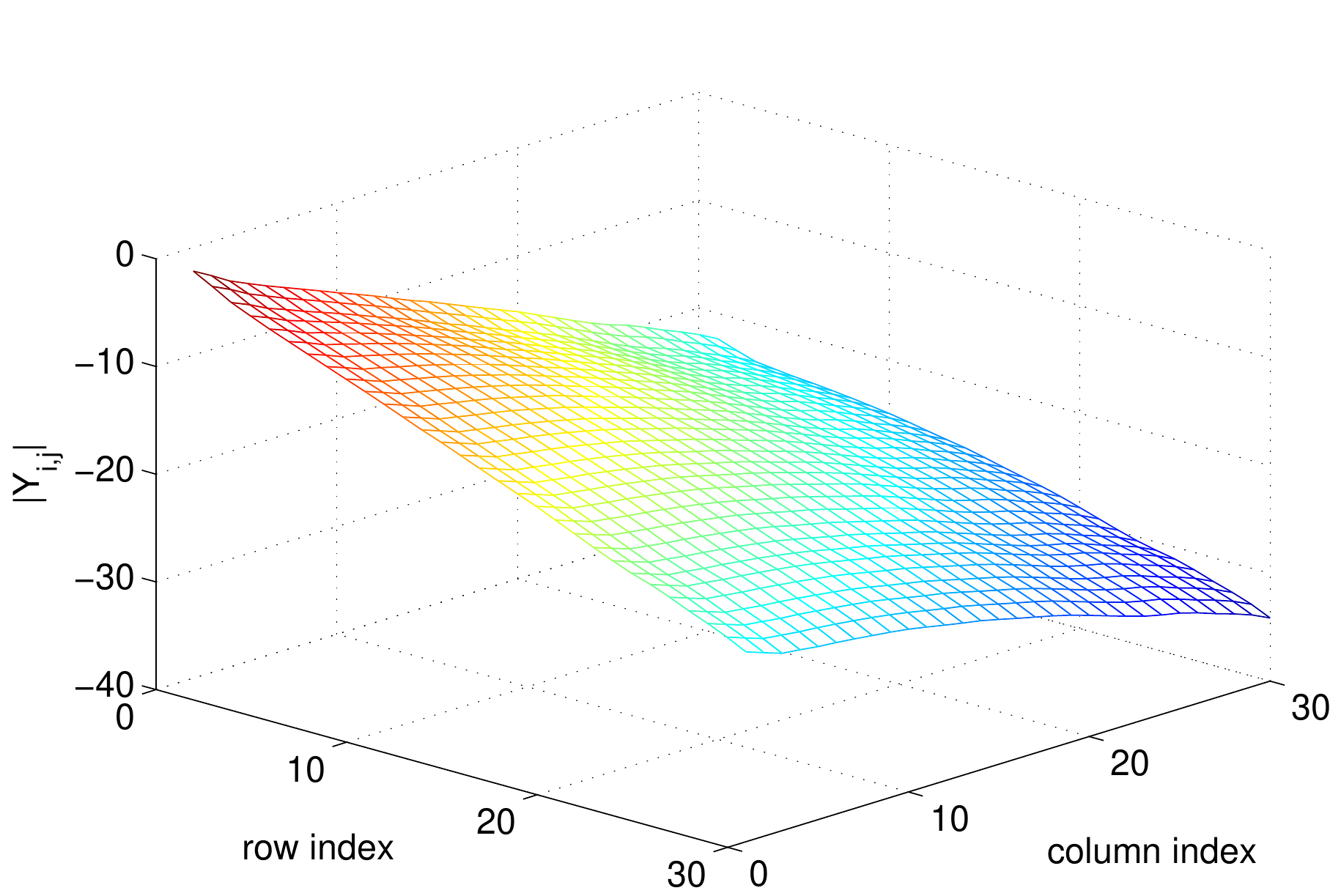}
\caption{Decay pattern of solution $Y$ of the reduced Lyapunov equation (\ref{eqn:reduced}). 
$A\in\RR^{n\times n}$, $n=900$ is the discretization of the
two-dimensional Laplacian in the unit square, $b$ is a random vector and $m=30$.
\label{fig:Y}}
\end{figure}

\begin{example}
{\rm
In Figure~\ref{fig:Y} a typical pattern is shown for $Y$: for this example, $A\in\RR^{n\times n}$,
$n=900$ is the
finite difference discretization of the two-dimensional Laplace operator in the unit square with
homogeneous boundary conditions, and $B=b$ is taken to be a vector with random uniformly distributed
values in the interval
$(0,1)$. A standard Krylov subspace of dimension $m=30$ was considered, so that $Y$
is a real symmetric $30\times 30$ matrix.
}
\end{example}

The analysis of the decay pattern of $Y$ appears to be new. On the other hand, it is
well known that the last row (or column) of $Y$ carries information on the accuracy of the
solution $X_m$. Indeed, let $\|R\|$ be the Frobenius ($\|\cdot \|_F$)
or induced 2-norm ($\|\cdot\|_2\|$) of the matrix $R$.
Then using (\ref{eqn:Arnoldi}) the residual norm satisfies (see \cite{Simoncini.survey13})
\begin{eqnarray*}
\|R\| &=& 
\| A V_m Y V_m^\TT + V_m Y V_m^\TT A^\TT - V_m e_1 \|b\|^2 e_1^\TT V_m^\TT \| \\
&=& 
\|  V_m T_m Y V_m^\TT + V_m Y T^\TT V_m^\TT  - V_m e_1 \|b\|^2 e_1^\TT V_m^\TT \\\
&& \qquad + v_{m+1} t_{m+1} e_m^\TT Y V_m^\TT + V_m Y e_m t_{m+1} v_{m+1}^\TT \| \\
&=& 
\|  V_{m+1} 
\begin{bmatrix}
 T_m Y  + V_m Y T^\TT   - e_1 \|b\|^2 e_1^\TT V_m^\TT &  t_{m+1} e_m^\TT Y\\
Y e_m t_{m+1} & 0 
\end{bmatrix} V_{m+1}^\TT \| \\
&=&
\|
\begin{bmatrix}
 0 &  t_{m+1} e_m^\TT Y\\
Y e_m t_{m+1} & 0 
\end{bmatrix}  \|,
\end{eqnarray*}
where in the last equality the orthogonality of the columns of $V_m$, and
(\ref{eqn:reduced}) were used. Hence,
$$
\|R\|_2 = |t_{m+1}| \, \|Y e_m \|_2, \quad
\|R\|_F =\sqrt{2}|t_{m+1}| \,  \|Y e_m \|_2 .
$$
Therefore, if the last column of $Y$, $Ye_m$, is small in norm then the residual will be small. We
note that on the other hand, the possibility
that $|t_{m+1}|$ is small is related to the fact that an invariant subspace is
found while generating the space $K_m$ via the Arnoldi recurrence 
(\ref{eqn:Arnoldi}), which is in general a less likely event, for general $b$ and
small $m$ compared to $n$.

Due to the previous discussion on the entries $|Y_{i,j}|$, for $j=m$ the vector $\|Ye_m\|$ 
is expected to be small for $m$ large enough (we recall here that $Y$ changes and
increases its dimension as $m$ increases). For $m>3$ the norm is bounded as
\begin{eqnarray*}
\|Ye_m\| & \le &
\|Ye_m\|_1 \\
& \le &
\frac{\|b\|^2}{2\pi}
\frac{64}{|\lambda_{\max}-\lambda_{\min}|^2} 
\sum_{i=1}^m\left( \int_{-\infty}^{\infty} 
 \left ( \frac{R^2}{(R^2-1)^2}\right )^2
 \left ( \frac{1}{R}\right )^{|i-1|+m-3} {\rm d}\omega \right) \\
& = &
\frac{\|b\|^2}{2\pi}
\frac{64}{|\lambda_{\max}-\lambda_{\min}|^2} 
 \int_{-\infty}^{\infty} 
\sum_{i=1}^m
 \left ( \frac{R^2}{(R^2-1)^2}\right )^2
 \left ( \frac{1}{R}\right )^{|i-1|+m-3} {\rm d}\omega  \\
& = &
\frac{\|b\|^2}{2\pi}
\frac{64}{|\lambda_{\max}-\lambda_{\min}|^2} 
 \int_{-\infty}^{\infty} 
 \left ( \frac{R^2}{(R^2-1)^2}\right )^2
 \left ( \frac{1}{R}\right )^{m-3} 
\sum_{k=0}^{m-1}
 \left ( \frac{1}{R}\right )^{k} {\rm d}\omega .
\end{eqnarray*}
Using $\displaystyle\sum_{k=0}^{m-1} \left(\frac 1 R\right )^k 
= (1-\frac 1 {R^m})/(1-\frac 1 R) = \frac{R}{R-1}\frac{ R^m-1}{R^m}$, after
some algebra we obtain
\begin{eqnarray*}
\|Ye_m\|  
&\le&
\frac{\|b\|^2}{2\pi}
\frac{64}{|\lambda_{\max}-\lambda_{\min}|^2} 
 \int_{-\infty}^{\infty} 
\frac{R^8}{(R-1)(R^2-1)^4}\frac 1 {R^m}  \frac {R^m-1} {R^{m}}  {\rm d}\omega \\
&<& 
\frac{\|b\|^2}{2\pi}
\frac{64}{|\lambda_{\max}-\lambda_{\min}|^2} 
 \int_{-\infty}^{\infty} 
\frac{R^8}{(R-1)(R^2-1)^4} \frac 1 {R^m} 
  {\rm d}\omega
\end{eqnarray*}
 where we used that $\frac {R^m-1} {R^{m}} <1$.
More explicit estimates could be obtained by further bouding the integrand
in terms of the variable $\omega$.

The relation between the convergence rate, in terms of the residual norm,
 and the solution sparsity pattern should not come as a surprise.
Indeed, the estimate in Theorem~\ref{th:freund} exploits the same
polynomial approximation properties 
that are used to prove convergence of standard Krylov subspace methods
to solve the Lyapunov equation by projection \cite{Simoncini.Druskin.09}.

\subsection{Numerical solution for a sparse right-hand side}
Projection type methods are suitable when the right-hand side matrix $D$ has
low rank, as discussed in the previous sections. Sparse but larger rank
right-hand sides can also enjoy particular properties, and may also be
viewed as sum of low rank matrices.
For instance, assume that $A$ is symmetric positive definite and
banded. Due to the linearity of the matrix equation,
the solution $X$ to
\begin{eqnarray}\label{eqn:Ddiag}
A X + X A = D, \quad D = {\rm diag}(\delta_1, \ldots, \delta_n), 
\end{eqnarray}
can be split as $X=X_1 + \cdots + X_n$, where 
for all $j\in\{1, \ldots n\}$ such that $\delta_j\ne 0$
the addend matrix $X_j$
 is the solution to the corresponding equation
$$
A X + X A =  e_j \delta_j e_j^\TT ;
$$
see, e.g., \cite{Hu.Reichel.92} for a similar linearity strategy.
It follows from our discussion on decay that each $X_j$ will have
an a-priori detectable decay pattern - a peak corresponding to the
$(j,j)$ entry -  from which the decay pattern of the whole of
$X$ can be predicted; see, e.g., \cite{Haber.Verhaegen.14tr} for
a similar discussion on the sparsity property of the resulting solution.
We remark that our decay analysis suggests that each $X_j$ may be
truncated to maintain sparsity in the approximate solution $\widetilde X_j$, so as to obtain
a good enough sparse approximate solution $\widetilde X$. We next report a typical 
example; see also Example~\ref{ex:trunc}.

\begin{figure}[htb]
\centering
\includegraphics[width=1.6in,height=1.6in]{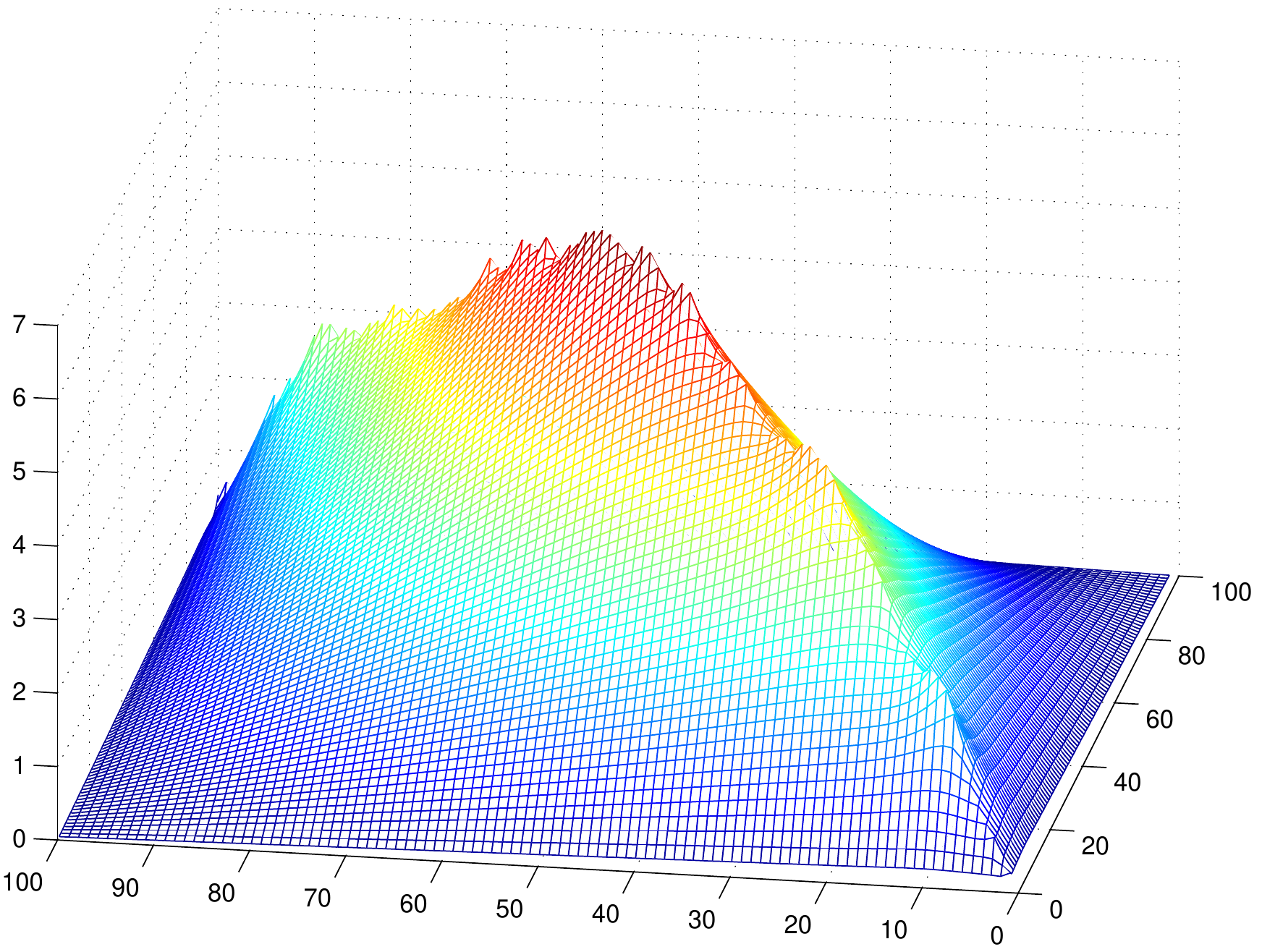}\hskip 0.5in
\includegraphics[width=1.6in,height=1.6in]{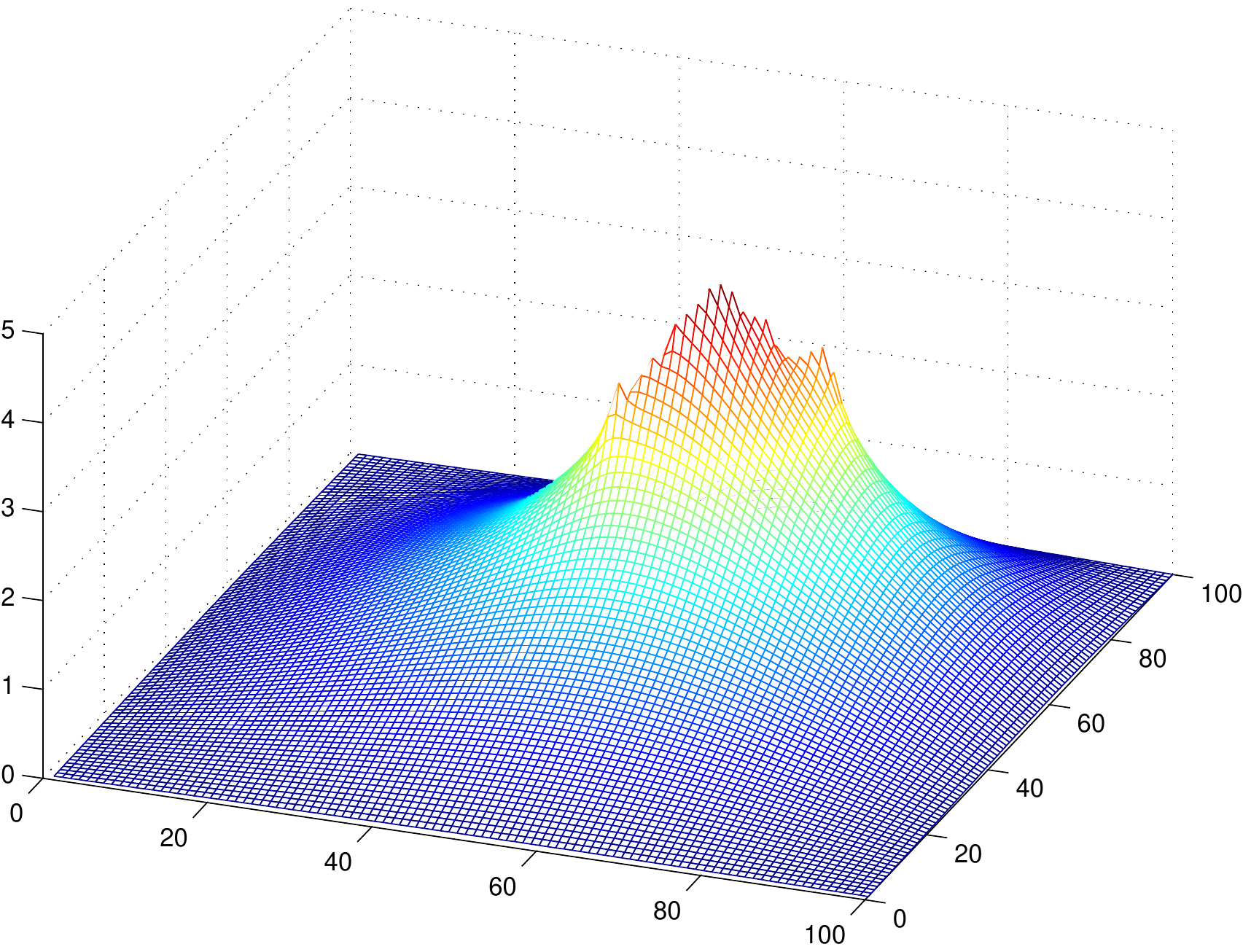}
\caption{Decay pattern of solution $X$ of the Lyapunov equation with 
right-hand side a diagonal matrix $D$. Left: $D$ nonsingular. Right: nonzero elements
in $D$ only corresponding to the diagonal entries from 50 to 70. 
A different viewpoint is used in the two plots.
\label{fig:Ddiag}}
\end{figure}

\begin{example}
{\rm
For $A$ equal to
the Laplace operator as before, the left plot of
Figure~\ref{fig:Ddiag} reports the solution to
(\ref{eqn:Ddiag}) when $D$ is a nonsingular diagonal matrix with
random entries uniformly distributed in the interval (0,1); the right
plot corresponds to $D={\rm diag}(\delta_1, \ldots, \delta_n)$ with $\delta_j$ equal
to a random value as above, but only for $j=50, \ldots, 70$.
The different sparsity patterns of the two solutions is as expected
from the theory. In particular, the left plot shows slightly larger values
of the diagonal entries corresponding to
the central part of the diagonal, where all linear components $X_j$ contribute.
On the other hand, the right plot confirms that (diagonal) peaks can only occur
in correspondence with the nonzero diagonal entries of $D$. 
}
\end{example}

\section{Further considerations}
Our presentation was aimed at highlighting some decay properties of the
matrices involved in the solution of the Lyapunov equation, that may give
insight into the development of new numerical methods. In particular,
the solution of the matrix equation when the right-hand side matrix $D$
has large or even full rank remains a big challenge.  Exploiting the possible
sparsity of $D$ may provide a good solution strategy; we refer the reader to
\cite{Haber.Verhaegen.14tr} and its references for recent developments in
this direction.

For ease of presentation,
our analysis was limited to the Lyapunov equation. 
Very similar results can be obtained for the Sylvester linear equation
$A_1X+XA_2 = D$, where $A_1, A_2$ may be of different dimensions; indeed, 
the solution can be written with a
closed form similar to that in (\ref{eqn:integral}), so that corresponding
decay bounds can be obtained; see \cite{Simoncini.survey13}.

Most of our results can be generalized to the case of $A$ nonnormal
and diagonalizable; this can be performed by replacing Theorem \ref{th:freund}
with corresponding results for functions of nonsymmetric matrices, developed
for instance in \cite{Benzi2007}.

\bibliography{%
/home/valeria/Bibl/Biblioteca}

\end{document}